\documentclass[a4paper,10pt]{article}
\usepackage[english]{babel}
\usepackage[T1]{fontenc} 
\usepackage{amsmath,amsfonts,amsthm,amssymb,bm} 
\usepackage{stmaryrd}

\usepackage{lipsum} 

\usepackage{sectsty} 

\usepackage{cite}
\usepackage[pdftex]{color,graphicx} 
\usepackage{color}    
\usepackage{here}
\usepackage{graphicx, subfigure}
\usepackage{enumerate}

\numberwithin{equation}{section} 
\numberwithin{figure}{section} 
\numberwithin{table}{section} 
\usepackage[usenames,dvipsnames,svgnames,table]{xcolor}

\usepackage[
    colorlinks=true, 
    citecolor=red,
    linkcolor=OliveGreen,
    urlcolor=blue]{hyperref}
    
\usepackage{proof}

\newtheorem{theorem2}{Theorem}
\newtheorem{theorem}[theorem2]{Theorem}

\newtheorem{lemma}[theorem2]{Lemma}
\newtheorem{claim}[theorem2]{Claim}

\newtheorem{prop}[theorem2]{Proposition}

\newtheorem{cor}[theorem2]{Corollary}

\numberwithin{theorem2}{section}

\theoremstyle{definition}
\newtheorem{defi}[theorem2]{Definition}
\newtheorem{remark}[theorem2]{Remark}
\newcommand\numberthis{\addtocounter{equation}{1}\tag{\theequation}}

\newcommand{\zfc}{{\rm ZFC}}

\newcommand{\Ord}{{\sf Ord}}
\newcommand{\Lim}{{\sf Lim}}

\newcommand{\lb}{\llbracket}
\newcommand{\rb}{\rrbracket}
\newcommand{\gl}{{\sf GL}}
\newcommand{\glp}{{\sf GLP}}
\newcommand{\bg}{{\sf BG}}

\newcommand{\csp}[1]{\tt cs(#1)}
\newcommand{\cs}{{\tt cs(\varsigma)}}
\newcommand{\otyp}{{\tt otyp}}

\newcommand{\polyt}{\{\mathcal{T}_\iota\}_{\iota < \Lambda}}

\newcommand{\ca}[1]{\mathcal{#1}}

\newcommand{\glpp}[1]{\glp_#1}
\newcommand{\glplam}{{\sf GLP_\Lambda}}

\newcommand{\langlam}{\mathcal L_\Lambda}

\newcommand{\ic}[2]{{#2}_{\uparrow {#1}}}
\newcommand{\ico}[1]{\mathcal{I}_{#1}}
\newcommand{\lift}{\Theta}

\usepackage{tikz}
\usetikzlibrary{graphs}
\usetikzlibrary{arrows}

\tikzset{
  treenode/.style = {align=center, inner sep=0pt, text centered,
    font=\sffamily},
  arn_n/.style = {treenode, circle, white, font=\sffamily\bfseries, draw=black,
    fill=black, text width=1.5em},
  arn_r/.style = {treenode, circle, red, draw=red, 
    text width=1.5em, very thick},
  arn_x/.style = {treenode, rectangle, draw=black,
    minimum width=0.5em, minimum height=0.5em}
}

\usepackage{authblk}
\title{A Topological Completeness Theorem for Transfinite Provability Logic\thanks{Email: {\texttt{aguilera@logic.at}}.
\newline
Mathematics Subject Classification 2010: 03F45, 03E10, 54G12. 
\newline 
Date: \today\, (compiled).}}
\author{J. P. Aguilera \\ 
Institute of Discrete Mathematics and Geometry,\\
Vienna University of Technology,\\
Wiedner Hauptstrasse 8-10,\\
1040 Vienna, Austria.}
\date{}

\begin{document}

\maketitle
\tableofcontents
\begin{abstract}
We prove a topological completeness theorem for the modal logic $\glp$ containing operators $\{\langle\xi\rangle:\xi\in\Ord\}$ intended to capture a wellordered sequence of consistency operators increasing in strength. More specifically, we prove that, given a tall-enough scattered space $X$, any sentence $\phi$ consistent with $\glp$ can be satisfied on a polytopological space based on finitely many Icard topologies constructed over $X$ and corresponding to the finitely many modalities that occur in $\phi$.
\end{abstract}

\clearpage
\section{Introduction}
The purpose of this article is to prove a topological completeness theorem for the transfinite extension of Japaridze's logic $\glp$. $\glp$ is a provability logic in a propositional language augmented with a possibly transfinite sequence of modal operators; our case of interest is that in which the sequence is wellfounded. Each of these operators can be interpreted arithmetically as asserting provability within a given theory, and the logic $\glp$ relates these notions of provability to one another. For arithmetical interpretations of $\glp$, see Beklemishev \cite{beklemishev04}, Fern\'andez-Duque and Joosten \cite{omegaruleinterpretation}, Cord\'on-Franco et al. \cite{cordonetal15}, and others. As a modal logic, $\glp$ has some unusual properties. For example, it is not complete with respect to any class of relational semantics; a natural question is whether it is complete with respect to its neighborhood (i.e., topological) semantics. Beklemishev and Gabelaia \cite{BekGab2014} showed that $\glp_\omega$, the restriction of $\glp$ to $\omega$-many modalities, is complete with respect to a natural topological space on the ordinal $\varepsilon_0$. Because sentences in the language of $\glp$ contain instances of only finitely many modalities, the spaces constructed by Beklemishev and Gabelaia serve as models also for formulas in the language of transfinite $\glp$; hence it is also topologically complete. However, it is an open problem whether transfinite $\glp$ is complete with respect to a single polytopological space equipped with a topology for each modality, but there are some very natural candidates. An example of these is what are known as the \emph{canonical topological semantics} for $\glp$. The question of completeness with respect to these spaces has very interesting connections with stationary reflection and indescribable cardinals (see Bagaria  \cite{bagariainf}, Bagaria \cite{Ba18}, Bagaria-Magidor-Sakai \cite{bagariaetal15} and Brickhill \cite{brickhillinf}). Completeness for the two-modality fragment with respect to these topologies was proved by Beklemishev \cite{completenessofglb}. It is not hard to see that $\glp$ is not strongly complete with respect to its canonical topological semantics; we shall prove this below.

Another example is the natural extension of the spaces from Beklemishev-Gabelaia \cite{topocompletenessofglp}. Completeness with respect to these spaces was proved by Fern\'andez-Duque \cite{polytopologies} for restrictions of $\glp$ to any countable amount of modalities. It is not known whether this result can be extended to arbitrarily long sequences of modalities, but the results from \cite{AguFer15} show that the techniques would need to be very different. 

The topological completeness theorem we shall prove here goes in this direction. Roughly, given a sentence consistent with $\glp$, we construct a topological model for it with finitely many topologies. The new feature is that these topologies in a way correspond to the modalities appearing in the sentence; we call these $\vec\vartheta$-polytopologies, for $\vec\vartheta$ a finite sequence of ordinals. In particular, one can extend the space with intermediate increasing topologies corresponding to modalities not appearing in the sentence in such a way that each topology results in a model of the unimodal $\gl$. Unfortunately, this extension (with the intermediate topologies) will not be a model of $\glp$, but we hope that a similar construction can yield models of $\glp$ and retain completeness. This hope is the main motivation for carrying out the work reported in this article.

Our main tool is a technical ``product lemma.'' Essentially, given two ordinals $\kappa$ and $\lambda$, we find an ordinal $\Theta$ and natural embeddings $\pi_0: \Theta\to\kappa$ and $\pi_1:\Theta \to\lambda$ which preserve satisfiability of polymodal formulas, if the ordinals are equipped with the right topologies. This is a generalization of a technical lemma of Beklemishev-Gabelaia \cite{topocompletenessofglp}, which corresponds to the case in which the first element of the sequence $\vec\vartheta$ is $1$. The proof is largely arithmetical and relies heavily on the theory of hyperexponentials and hyperlogarithms of Fern\'andez-Duque and Joosten \cite{hyperations}.

\section{Preliminaries}\label{TopModels}
\subsection{The polymodal logic of provability}
For any ordinal number $\Lambda$ we consider a language $\langlam$ consisting of a countable set of propositional variables $\mathbb{P}$ together with the constants $\top$, $\bot$; Boolean connectives $\land$, $\lor$, $\lnot$, $\rightarrow$; and a modality $[\xi]$ for each ordinal $\xi < \Lambda$. As usual, we write $\langle \xi \rangle$ as a shorthand for $\lnot [\xi] \lnot$.  
\begin{defi}
The logic $\glplam$ is then defined to be the least logic containing all propositional tautologies and the following axiom schemata: 

\begin{enumerate}
\item[(i)] $[\xi](\varphi \rightarrow \psi)\rightarrow([\xi]\varphi \rightarrow [\xi]\psi)$ for all $\xi < \Lambda$,
\item[(ii)] $[\xi]([\xi]\varphi \rightarrow\varphi)\rightarrow[\xi]\varphi$ for all $\xi < \Lambda$,
\item[(iii)] $[\xi]\varphi \rightarrow [\zeta]\varphi$ for all $\xi < \zeta < \Lambda$,
\item[(iv)] $\langle \xi \rangle \varphi \rightarrow [\zeta] \langle \xi \rangle \varphi$ for all $\xi < \zeta < \Lambda$,
\end{enumerate}
and closed under the rules modus ponens and necessitation for each $[\xi]$:
\begin{align*}
\infer[\mathrm{MP}]{\psi}{\varphi \rightarrow \psi \qquad \varphi}&  &
\infer[\mathrm{nec}] {[\xi] \varphi}{\varphi}
\end{align*}
\end{defi}

We will often write simply $\glp$ for $\glplam$ when we do not want to specify a $\Lambda$.
Note that $\glp$ when restricted to any one modality is simply the well-known logic $\gl$. Modal logics are usually studied by means of relational semantics. A \textit{Kripke $\Lambda$-frame} is a structure $(W, \{R_{\xi}\})_{\xi < \Lambda}$, where each $R_\xi$ is a binary relation on $W$. We define a \textit{valuation} $\lb \cdot \rb$ to be a function assigning subsets of $W$ to each $\langlam$-formula such that $\lb \cdot \rb$ respects boolean connectives and such that
\begin{equation*}
\lb \langle \xi \rangle \varphi \rb = R_\xi^{-1} \lb \varphi \rb.
\end{equation*}

\begin{prop}[Segerberg \cite{segerberg1971}] \label{segerbergs theorem}
$\glpp 1$ is complete with respect to the class of finite relational structures $(W, R)$ that are conversely wellfounded trees.
\end{prop}

The preceding proposition provides a convenient way to study $\glpp 1$. However, as is well known, $\glplam$ is incomplete with respect to any class of relational structures whenever $1 < \Lambda$. This motivates the search for topological models of the $\glp$.

Recall that $x$ is a \textit{limit point} of $A$ if $A$ intersects every punctured neighborhood of $x$. We call the set of limit points of $A$ the \textit{derived set} of $A$ and denote it by $dA$. We may also denote it by $d_\tau A$ to emphasize the topology we are considering. The derived set operator is iterated transfinitely by setting
\begin{enumerate}
\item $d^0A=A$, 
\item $d^{\alpha + 1}A = d d^\alpha A$, and 
\item $d^{\gamma}A = \bigcap_{\alpha < \gamma}d^\alpha A$ for limit ordinals. 
\end{enumerate} 
Since $d^\alpha X \supset d^\beta X$ whenever $\alpha < \beta$, there exists a minimal ordinal $ht(X)$---the \textit{height} or \emph{rank} of $X$---such that $d^{ht(X)}X = d^{ht(X) + 1}X$. For any $x \in X$, we let $\rho_\tau x$, the \textit{rank} of $x$, be the least ordinal $\xi$ such that $x \not \in d^{\xi+1} X$, if it exists. 

Throughout this paper, we will speak about \textit{rank-preserving} extensions of topologies. 

\begin{lemma}[Beklemishev-Gabelaia \cite{topocompletenessofglp}]\label{CharRPE}
A topology $\sigma$ is a rank-preserving extension of a scattered topology $\tau$ if, and only if, $\rho_\tau[U]$ is an ordinal for each $U \in \sigma$.
\end{lemma}

A point in $A$ that is not a limit point is \textit{isolated}. Thus a point is isolated if and only if it has rank $0$. We denote by $iso(A)$ the set of isolated points in $A$. A topological space is \textit{scattered} if $iso(A) \neq \varnothing$ for each $A \subset X$ (alternatively, if $d^{ht(X)}X = \varnothing$). Not all scattered spaces are $T_1$ (e.g., $X = \{0,1\}$ with open sets $\varnothing$, $\{0\}$ and $X$), however, the examples in which we will focus are.

We study \textit{polytopological spaces}---structures $(X, \polyt)$, where $X$ is a set and $\polyt$ is a sequence of topologies of length $\Lambda$. Topological semantics for modal logics may be defined by interpreting diamonds as topological derivatives.

\begin{defi}[Topological semantics]\label{DefKripkSem}
Let $\mathfrak{X}=(X, \polyt)$ be an polytopological space. A \textit{valuation} is a function $\lb\cdot\rb: \langlam \to \mathcal{P}(X)$ such that for any $\langlam$-formulae $\varphi, \psi$:
\begin{enumerate}
\item[(i)] $\lb\bot\rb = \varnothing$;
\item[(ii)] $\lb\lnot \varphi\rb = X \backslash \lb\varphi\rb$;
\item[(iii)] $\lb\varphi \land \psi\rb = \lb\varphi\rb \cap \lb\psi\rb$;
\item[(iv)] $\lb \langle \xi \rangle \varphi\rb = d_{\mathcal{T}_\xi} \lb\varphi\rb$.
\end{enumerate}

A \textit{model} $\mathfrak{M} = (\mathfrak{X}, \lb\cdot\rb)$ is a polytopological space together with a valuation.  We say that $\varphi$ is \textit{satisfied} in $\mathfrak{M}$ if $\lb\varphi\rb$ is nonempty and we say $\varphi$ is \textit{valid} in a space $\mathfrak{X}$ and write $\mathfrak{X} \models \varphi$ if $\lb\varphi\rb = X$ for any model based on $\mathfrak{X}$. 
\end{defi}

In order that a space validate the axioms of $\glp$, we need to impose some regularity conditions (see Beklemishev-Bezhanishvili-Icard \cite{BekBezIc}). A space $(X, \polyt)$ is a \textit{$\glplam$-space} if $\polyt$ is non-decreasing, scattered, and 

\begin{equation}\label{CondAmb}
d_\xi A \in \ca{T}_\zeta \text{ for all } \xi < \zeta \text{ and all } A \in \mathcal{P}(X)
\end{equation}
In the situation above, we refer to $\polyt$ as a \textit{$\glplam$-polytopology}. Clearly, we have:

\begin{lemma}
Any $\glplam$-space validates all theorems of $\glplam$.
\end{lemma}

A natural way of constructing $\glplam$-polytopologies appears to be to start with any scattered topology and simply add all derived sets at each stage, thus making $\ca A = \wp(X)$. This results in what has come to be known as the \textit{canonical $\glp$-space generated by $X$}. In doing so, the topologies quickly become extremely fine. In fact, for the most natural examples, their non-discreteness becomes undecidable within $\zfc$ after two or three iterations.

One way out of this, explored in Beklemishev-Gabelaia\cite{topocompletenessofglp}, is to extend the topology at each stage before adding derived sets. Extending the topology reduces the amount of derived sets attainable and makes subsequent topologies coarser. A different approach, introduced in Fern\'andez-Duque \cite{polytopologies}, is to fix increasing topologies from the beginning and restrict the algebra of possible valuations. We will consider the first approach here.

Let us finish this section with the remark that $\glp$ is not strongly complete with respect to its canonical semantics.
By \textit{strong completeness} (with respect to a class of models $\mathcal X$), we mean the following assertion: whenever $\Gamma$ is a set of $\langlam$-sentences consistent with $\glplam$, then there is some model $\mathfrak X \in \mathcal X$ where $\Gamma$ is satisfied. 
\begin{prop}\label{PropGLPNotStronglyComplete}
Suppose $X$ is a scattered $T_1$ space in which every $G_\delta$ set is open. Then $\gl$ is not strongly complete with respect to $X$.
\end{prop}
\proof
This is a generalization of the usual proof that $\gl$ is not strongly complete with respect to trees. Let 
$$\Gamma = \{\Diamond p_0\} \cup \{\Box(p_i \to \Diamond p_{i+1}) : i < \omega\}.$$
Suppose $\Gamma$ is satisfied at some $x$. Then, for each $i$, there is a punctured neighborhood $U_i$ of $x$ in which any point satisfying $p_i$ is a limit of points satisfying $p_{i+1}$. Since $X$ is $T_1$, each $U_i$ is open. It follows that $U = \bigcap_{i < \omega} U_i$ is $G_\delta$ and thus open. Since $x$ satisfies $\Diamond p_0$, $U$ contains some $x_0$ of some rank $\alpha_0$ satisfying $p_0$. Inductively, for each $i < \omega$, there is some $x_i \in U$ satisfying $p_i$ and, by 
$$\Box(p_i \to \Diamond p_{i+1}),$$ 
$U$ contains some $x_{i+1}$ of rank $\alpha_{i+1} < \alpha_i$ satisfying $p_{i+1}$. This gives an infinite decreasing sequence of ordinals.
\endproof

Recall that if $\kappa$ is an ordinal of uncountable cofinality, then the intersection of countably many sets which are closed and cofinal in $\kappa$ is also closed and cofinal in $\kappa$. Hence, 
Proposition \ref{PropGLPNotStronglyComplete} applies to the closed-unbounded topology from Blass \cite{blass1990}. More generally:

\begin{cor}
$\gl$ is not strongly complete with respect to topologies on ordinals given by countably complete filters, such as the closed-unbounded topology.
\end{cor}

The spaces we will consider will instead be built around the \emph{Generalized Icard topologies}.
\begin{defi}[Generalized Icard Topologies] \label{DefGIT}
Let $(X, \tau)$ be a scattered space of rank $\Theta$. We define a topology $\ic 1\tau$ generated by $\tau$ and all sets of the form
\begin{equation*}
(\alpha, \beta)^{\tau} := \{x \in X: \alpha < \rho_\tau x < \beta\},
\end{equation*}
for ordinals $\alpha < \beta \leq \Theta + 1$. We iterate this construction by setting
\begin{itemize}
\item $\ic{(\iota + 1)}\tau = \ic{1}{(\ic{\iota}\tau)}$, and
\item $\ic\lambda\tau = \bigcup_{\xi < \lambda}\ic\xi\tau$ at limit stages.
\end{itemize}
These are called the \textit{generalized Icard topologies}.
\end{defi}

These topologies were defined differently in \cite{AguFer15}. By Lemma \ref{properties of icard spaces}.\ref{CharIcard} below, both definitions coincide. Another equivalent formulation is as follows: $\ic 1 \tau$ is generated by $\tau$ and the family
\begin{equation} \label{AltDefGIT}
\{d^\xi X\colon \xi \in \Ord\}.
\end{equation}

\subsection{Arithmetic, I}
\begin{defi}
We fix some notation related to ordinal arithmetic.
\begin{enumerate}
    \item Whenever $\alpha < \beta$, we denote by $-\alpha + \beta$ the unique ordinal $\gamma$ such that $\alpha + \gamma = \beta$. 
    
    \item Whenever $A$ is a set of ordinals, we denote by $\alpha + A$ the set $\{\alpha + \beta \colon \beta \in A\}$. Expressions such as $-\alpha + A$ are defined analogously, if they make sense.
        
    \item For all nonzero $\xi$, there exist ordinals $\alpha$ and $\beta$ such that $\xi = \alpha + \omega^\beta$. Such a $\beta$ is unique. We denote it by $\ell\beta$ and call it the \emph{end-logarithm} of $\xi$.
            
    \item For all nonzero $\xi$, there exists a unique ordinal $\eta$ such that $\xi$ can be written as $\omega^\eta + \gamma$, with $\gamma < \xi$. We denote this ordinal by $L\xi$ and call it the \emph{initial logarithm} of $\xi$.
    
\end{enumerate}
\end{defi}

The operations $\ell$ and $L$ should be regarded as functions on (a sufficiently large subset of) $\Ord$. Nonetheless, in its use and in general whenever we deem it convenient, we will omit the symbol `$\circ$' $ $ for function composition, as well as perhaps parentheses.

Our completeness proof will rely heavily upon an analysis of generalized Icard topologies and their structure induced by the arithmetical properties of ordinals. Hence, developing a thorough intuition about them will be crucial. A most useful remark in this direction is the fact that they are to arbitrary topological spaces as the usual order topology is to ordinal numbers. Indeed, define the \textit{initial segment topology} $\ico 0$ on an ordinal $\Theta$ (or on $\Ord$) to be generated by all initial segments $[0, \alpha]$, for $\alpha < \Theta$. Then $(\Theta, \ico 0)$ is a scattered space: a rather trivial scattered space---it carries no further information than the usual ordering on $\Ord$. For instance, we have $\rho_{\ico 0}\alpha = \alpha$ for all $\alpha$ and $ht(\Theta, \ico 0) = \Theta$.

\begin{lemma} \label{LemmIntTopo}
$\ico 1 := \ic{1}{\ico 0}$ is the order topology. We have $\rho_{\ico 1}\alpha = \ell\alpha$ for all $\alpha$, so in particular isolated points are exactly the successor ordinals. Moreover, $ht(\Theta, \ico 1) = L(\Theta) + 1$.
\end{lemma}
\proof
It is not hard to see that $\ico 1$ is the order topology, and that the rank function is $\ell$ is established by a simple induction. Finally, let $\mathbb H$ be the class of additively indecomposable ordinals. It follows that
\begin{align*}
ht(\Theta, \ico 1) 
= \sup_{\xi < \Theta}(\ell\xi + 1)
= \sup_{\xi < \Theta \cap \mathbb{H}}(\ell\xi + 1) 
= \sup_{\xi < \Theta \cap \mathbb{H}}(L\xi + 1) 
= L(\Theta) + 1,
\end{align*}
as claimed.
\endproof

In what follows, we write simply $\ico \lambda$ instead of $\ic \lambda{\ico 0}$. 
These topologies are important because,
as we will see, the completeness theorem can quickly be reduced to the case when the underlying space is an ordinal equipped with a topology of the form $\ico\lambda$.

\subsection{d-maps and J-maps}
There is an appropriate notion of structure-preserving mappings between scattered spaces. We say that a function between topological spaces is \textit{pointwise discrete} if the preimage of any singleton is a discrete subspace.

\begin{defi}[d-map]
Let $X$ and $Y$ be scattered spaces. A function $f:X \to Y$ is a \textit{d-map} if it is continuous, open, and pointwise discrete.
\end{defi}

Clearly, any homeomorphism is a d-map. In particular, ordinal addition and substraction, i.e., functions of the form
\begin{equation*}
(-\xi + \cdot)\colon ([\xi,\xi+\Theta], \ico\zeta) \to ([0,\Theta], \ico\zeta),
\end{equation*}
are d-maps. The rank function 
$$\rho_\tau\colon(X,\tau) \to ([0,\rho_\tau X],\ico 0)$$
is also a d-map. A more interesting example is given by end-logarithms of the form:
\begin{equation*}
\ell \colon (\Theta, \ico{1 +\zeta}) \to (\Theta, \ico\zeta).
\end{equation*}
A proof of this, and the more general Lemma \ref{logdmaps} below can be found in Fern\'andez-Duque \cite{polytopologies}.

Since the composition of d-maps is a d-map, they can be thought of as morphisms in the category of scattered spaces. We will now state various properties of d-maps. 

\begin{lemma} \label{propertiesofdmaps}
Let $f:X \to Y$ be a d-map.
\begin{enumerate}
\item \label{uniquedmap} If $Y$ is an ordinal $\Theta$ with the initial segment topology, then $f$ is the rank function on $X$.
\item \label{PropDmaps2} For any $A \subset Y$, $f^{-1}dA = df^{-1}A$.
\item \label{nextlambdamap} $f:(X, \ic\lambda\tau) \to (Y, \ic\lambda\sigma)$ is a d-map for any $\lambda$.
\item \label{PropDmaps4} If $f$ is surjective, then for any $\langlam$-formula $\varphi$, $X \models \varphi$ implies $Y \models \varphi$.
\end{enumerate}
\end{lemma}
\proof
Items \ref{uniquedmap} and \ref{PropDmaps2} appear in Beklemishev-Gabelaia \cite{topocompletenessofglp}; item \ref{PropDmaps4} appears in Bezhanishvili-Mines-Morandi \cite{bmm} in the current formulation. Item \ref{nextlambdamap} is proved in \cite{AguFer15}, but therein a different definition of $\ic\lambda\tau$ is used, and we still have not shown that they are equivalent. Nonetheless, the claim can be proved by an easy induction. 

Suppose $f:(X, \ic\lambda\tau) \to (Y, \ic\lambda\sigma)$ is a d-map. Note that \ref{PropDmaps2} implies that d-maps are rank-preserving, i.e., 
\begin{equation*}
\rho_{\ic\lambda\tau}(x) = \rho_{\ic\lambda\sigma}(f(x)) \text{ for every } x \in X.
\end{equation*}
It follows that for any $\ic\lambda\tau$-open $A$, 
\begin{equation*}
f(A \cap (\alpha,\beta)^{\ic\lambda\tau}) = f(A) \cap (\alpha,\beta)^{\ic\lambda\sigma} \in \ic{\lambda+1}\sigma;
\end{equation*}
and for any $\ic\lambda\sigma$-open $B$, 
\begin{equation*}
f^{-1}(B \cap (\alpha,\beta)^{\ic\lambda\sigma}) = f^{-1}(B) \cap (\alpha,\beta)^{\ic\lambda\tau} \in \ic{\lambda+1}\tau;
\end{equation*}
so that $f$ is $(\ic{\lambda+1}\tau, \ic{\lambda+1}\sigma)$-continuous and open. Clearly it is also pointwise discrete. The case for limit $\lambda$ follows from Fern\'andez-Duque \cite[Lemma 5.8]{polytopologies}.
\endproof

As mentioned in the proof of \ref{propertiesofdmaps}.\ref{nextlambdamap}, Lemma \ref{propertiesofdmaps}.\ref{PropDmaps2} implies that d-maps are rank-preserving. Also, it follows from \ref{propertiesofdmaps}.\ref{nextlambdamap} that if the rank of $(X, \tau)$ is $\Theta$, then 
$$\rho_\tau: (X, \ic\lambda\tau) \to (\Theta, \ico\lambda)$$
is a d-map. The main feature of d-maps is as follows: 

\begin{lemma} \label{Lemmdmapsiff}
$\glp_1$ is complete with respect to a scattered space $(X, \tau)$ if, and only if, for any finite, converse-wellfounded tree $T$, there exists a $\tau$-open subspace $S$ of $X$ and a d-map $f\colon (S, \tau) \to T$.
\end{lemma}
\proof
That completeness follows from the existence of d-maps is independently due to Abashidze \cite{abashidze1985} and Blass \cite{blass1990}. Note that it immediately follows from Proposition \ref{segerbergs theorem} and Lemma \ref{propertiesofdmaps}.\ref{PropDmaps4}. 

The converse is probably folklore and will not be needed below, but we prove it anyway. Suppose  $\glp_1$ is complete with respect to $(X, \tau)$, where $\tau$ is scattered. Let $(T, <)$ be a finite, converse wellfounded tree. We define from $T$ a modal formula $\varphi$ consistent with $\glp_1$. Let $\{p_t\colon t\in T\}$ be a set of distinct propositional variables and $r$ be the root of $T$. Set
\begin{align*}
\varphi 
=\ & p_r 
\land \bigwedge_{s \in T;\, s \neq r}\neg p_s 
\land \left( \bigwedge_{t \in T; t \neq r}\Diamond p_t \right)
\land \Box\left( \bigvee_{t \in T}p_t\right) 
\land \Box\left( \neg p_r \right) \\
&\land \left( \bigwedge_{s,t \in T;\, s \neq t} \Box(p_s \rightarrow \neg p_t) \right) 
\land \left( \bigwedge_{s < t} \Box(p_s \rightarrow \Diamond p_t)\right) \\
&\land \left( \bigwedge_{s \not < t} \Box(p_s \rightarrow \neg \Diamond p_t) \right) 
\land \Box\bigwedge_{t\in T} \left(  p_t \rightarrow \Box \bigvee_{t < s}p_s\right).
\end{align*}
Clearly, there is a Kripke model based on $T$ where $\varphi$ is true in $r$; namely, any one where each $p_t$ holds only in $t$. Hence, $\varphi$ is consistent with $\glp_1$, whereby it is satisfiable in $X$. Fix a valuation over $X$ and a point $x_r \in X$ such that $x_r \models \varphi$. Thus, $x_r$ satisfies $p_r$ and $\bigwedge_{s,t \in T;\, s \neq r}\neg p_s$ and $x_r$ is a limit point of points satisfying each of $p_t$, for $t \neq r$. Moreover, by each of the conjuncts above:
\begin{enumerate}[i.]
\item there is a punctured neighborhood of $x_r$ where each point satisfies $p_t$ for some $t \in T$;
\item there is a punctured neighborhood of $x_r$ where no point satisfies $p_r$;
\item for each pair of distinct $s,t \in T$, there is a punctured neighborhood of $x_r$ of points satisfying at most one of $p_s$ and $p_t$;
\item for each pair of distinct $s,t \in T$ with $s<t$, there is a punctured neighborhood of $x_r$ where all points satisfying $p_s$ are limits of points satisfying $p_t$;
\item for each pair of distinct $s,t \in T$ with $s\not<t$, there is a punctured neighborhood of $x_r$ where all points satisfying $p_s$ are not limits of points satisfying $p_t$; and
\item there is a punctured neighborhood of $x_r$ where whenever a point $x$ satisfies $p_t$, then there is a punctured neighborhood of $x$ where each points satisfies one of $p_s$, with $t < s$.
\end{enumerate}
Let $S$ be the intersection of all those finitely many open neighborhoods of $x_r$. Clearly, $x \models p_t$ and $t \neq s$ together imply $x \not \models p_s$. We define $f\colon S \to T$ by
\begin{equation*}
f(x) = t \text{ if, and only if, } x\models p_t.
\end{equation*}

We claim $f$ is a d-map. Let $A_t$ be an open subset of $T$ of the form 
$$\{s \in T\colon t \leq s\},$$ 
so that 
$$f^{-1}(A_t) = \{x \in S\colon x \models p_s \land t \leq s\}.$$
This clearly equals $S$ if $t = r$. Otherwise, for each $x \in S$ with $x \models p_s$ and $t \leq s$, there is an open neighborhood $U$ of $x$ where each point satisfies $p_u$ for some $u > s$. But $t < u$, whence $U \subset f^{-1}(A_t)$. This implies that $f^{-1}(A_t)$ is open, and so $f$ is continuous.

Conversely, suppose $U \subset S$ is open, $x \in U$ is such that $x \models p_t$, and $t < s$. Then $x$ is a limit of points satisfying $p_s$, so that 
$$\{y \in S\colon y\models p_s\} \cap U \neq \varnothing,$$
whence $s \in f(U)$. Hence, $f$ is open. Finally if $t \in T$, then $f^{-1}(t)$ is discrete, for $t$ is the image of points satisfying $p_t$ and no point in $S$ can satisfy $p_s \land \Diamond p_s$ for any $s$. Therefore, $f$ is a d-map.
\endproof

Hence, the need to check whether a given space $\mathfrak X$ satisfies a formula is replaced by the definition of a suitable mapping between $\mathfrak X$ and some other space which is known to do so. In practice, polymodal analogs of Lemma \ref{Lemmdmapsiff} do not even require us to use full d-maps, but rather a weaker form of embeddings, as shown by Beklemishev \cite{kripkesemanticsofglp}:

\begin{defi}[J-frame]
A finite polymodal Kripke frame 
$$\mathfrak{F}=(W, \{<_n\}_{n<\omega})$$ 
is called a \textit{J-frame} if each relation is transitive and conversely wellfounded and it satisfies the following two conditions:
\begin{enumerate}
\item[(I)] For all $x, y \in W$ and all $m < n$: $x <_n y$ implies that for all $z \in W$: $x <_m z$ if, and only if $y <_m z$.
\item[(J)] For all $x, y, z \in W$ and all $m < n$: if $x <_m y$ and $y <_n z$, then $x <_m z$.
\end{enumerate}

We call a J-frame a \textit{J$_n$-frame} if all binary relations past the $n$th one are empty.
\end{defi}

Let $(T, <_0, \hdots, <_N)$ be a frame. Denote by $E_n$ the reflexive, symmetric, and transitive closure of $\bigcup_{n \leq k < \omega}<_k$. The equivalence classes under $E_n$ are called \textit{$n$-planes}. A natural order is defined on the set of $(n+1)$-planes:
\begin{equation*}
\alpha \prec \beta \text{ if, and only if, } x<_ny \text{ for some } x \in \alpha, y \in \beta
\end{equation*}

We say that a J-frame is a \textit{J-tree} if for all $n$, the $(n+1)$-planes contained in each $n$-plane form a tree under $<_n$ and if whenever $\alpha < \beta$ for two $(n+1)$-planes $\alpha, \beta$, we have $x <_n y$ for all $x\in \alpha$ and $y \in \beta$. This means that each J$_n$-tree can be thought of as a tree each of whose nodes is itself a J$_{n-1}$-tree. Below, a node $t \in T$ is a \emph{hereditary $k$-root} if for no $j\geq k$ and no $s \in T$ do we have $s<_j t$. We also write $x\ll_k y$ if $x<_j y$ for some $j \geq k$ and
\[{\ll_k}(x) = \{x \in T: x\ll_k y\}.\]

\begin{defi}[J-map]
Let $(T, \sigma_0, ... ,\sigma_n)$ be a J$_n$-tree and $(X, \tau_0, ... ,\tau_n)$ be a space with $n+1$ topologies. We say that a function $f:X \to T$ is a \textit{J$_n$-map} if
\begin{enumerate}
\item[($j_1$)] $f:(X,\tau_n) \to (T,\sigma_n)$ is a d-map;
\item[($j_2$)] $f:(X,\tau_k) \to (T,\sigma_k)$ is open for each $k$;
\item[($j_3$)] $f^{-1}({\ll_k}(x)),f^{-1}(\{x\}\cup {\ll_k}(x)) \in \tau_k$ for each $k<n$ and each hereditary $(k+1)$-root $x$;
\item[($j_4$)] $f^{-1}x$ is a $\tau_k$-discrete subspace for each $k<n$ and each hereditary $(k+1)$-root $x$.
\end{enumerate}
\end{defi}

\begin{lemma}[Beklemishev \cite{kripkesemanticsofglp}]\label{jmapcomposition}
If $f: Y \to Z$ is a J$_n$-map and $g: X \to Y$ is a d-map, then $f \circ g: X \to Z$ is a J$_n$-map.
\end{lemma}

\begin{lemma}[Beklemishev-Gabelaia \cite{topocompletenessofglp}]\label{to comp 4}
For each $\mathcal{L}_n$-formula $\varphi$ consistent with $\glp$, there exists a J$_n$-tree $T$ such that if $\mathfrak X$ is a $\glp_n$-space and $f: \mathfrak X \to T$ is a surjective J$_n$-map, then $\mathfrak{X} \models \varphi$. 
\end{lemma}

We call the tree obtained in Lemma \ref{to comp 4} the \textit{canonical tree for $\varphi$}.

\subsection{Arithmetic, II}
We will need the definition of hyperlogarithms and hyperexponentials, due to Fern\'andez-Duque and Joosten  \cite{hyperations}:
\begin{defi}\label{ExpDef}\
\begin{enumerate}
    \item \label{ExpDefPart2} The \textit{hyperlogarithms} $\{\ell^\xi\}_{\xi \in \Ord}$ are the unique family of pointwise maximal initial\footnote{That is, sending initial segments of $\Ord$ onto initial segments of $\Ord$.} functions that satisfy:
    \begin{enumerate}
    \item $\ell^1 = \ell$, and
    \item \label{LogDefb}$\ell^{\alpha + \beta} = \ell^\beta  \ell^\alpha$.
    \end{enumerate}

    \item \label{ExpDefPart1} Let the function $e$ be defined by $\xi \mapsto -1 + \omega^\xi$. The          \textit{hyperexponentials} $\{e^\zeta\}_{\zeta \in \Ord}$ are the unique pointwise minimal family of normal functions that satisfy
        \begin{enumerate}
        \item $e^1 = e$, and
        \item\label{ExpDef2} $e^{\alpha + \beta} = e^\alpha e^\beta$ for all $\alpha$ and $\beta$.
        \end{enumerate}
\end{enumerate}
\end{defi}

\noindent One can verify by induction that the sequence $\{\ell^\xi\gamma\colon \xi \in \Ord\}$ is non-increasing for any ordinal $\gamma$. If we set $e^0$ to be the identity function and $e^\xi 0 = 0$ for all $\xi$, then one can also describe hyperexponentials recursively by condition \ref{ExpDef}.\ref{ExpDef2}, together with the following normality clause:
\begin{flalign} \label{normality}
\text{for any } \xi \text{ and any limit }\lambda: e^\xi\lambda = \lim_{\eta \to \lambda}e^\xi\eta;&&
\end{flalign}
and the following fixed-point clause:
\begin{flalign} \label{FixPoints}
\text{for any } \xi \text{ and any limit }\lambda: e^\lambda(\xi+1) = \lim_{\eta \to \lambda}e^\eta(e^\lambda\xi + 1).&&
\end{flalign}

The hyperexponential family refines the Veblen hierarchy. We mention some more properties of hyperlogarithms and --exponentials.

\begin{lemma}[see Fern\'andez-Duque and Joosten \cite{hyperations} and Fern\'andez-Duque \cite{polytopologies}]\ \label{properties of exponentials}
\begin{enumerate}
\item If $\xi$ and $\delta$ are nonzero, then $\ell^\xi(\gamma + \delta) = \ell^\xi\delta$; if $\gamma < \delta$ as well, then $\ell^\xi(-\gamma + \delta) = \ell^\xi\delta$. Moreover, if $1 < \xi$, then $\ell^\xi(\gamma\delta) = \ell^\xi\delta$; \label{SumCancel}

\item If $\xi < \zeta$, then $\ell^\xi e^\zeta = e^{-\xi + \zeta}$ and $\ell^\zeta e^\xi = \ell^{-\xi + \zeta}$. Furthermore, if $\alpha < e^\xi\beta$, then $\ell^\xi \alpha < \beta$. \label{exponentialcancellation}
\end{enumerate}
\end{lemma}

\proof[Sketch of \ref{SumCancel}]
That $\ell^\xi(\gamma + \delta) = \ell^\xi\delta$ is proved by induction on $\xi$ using \ref{ExpDef}.\ref{LogDefb}. From this follows that if $\gamma < \delta$, then
$$\ell^\xi(\delta) = \ell^\xi(\gamma + (-\gamma + \delta)) = \ell^\xi(-\gamma + \delta).$$ 
Finally, it is proved by induction that $\ell(\gamma\delta) = L\gamma + \ell\delta$, so that if $1 < \xi$, then 
\begin{equation*}
\ell^\xi(\gamma\delta) 
= \ell^{-1 + \xi}\ell(\gamma\delta) 
= \ell^{-1 + \xi}(L\gamma + \ell\delta)
= \ell^{-1 + \xi}\ell\delta
= \ell^\xi\delta,
\end{equation*}
as desired.
\endproof

We now give an alternative characterization of topologies $\ic\lambda\tau$ and their rank functions:

\begin{lemma}\label{properties of icard spaces}
Let $(X, \tau)$ be a scattered space of rank $\Theta$.
\begin{enumerate}
\item \label{CharIcard} The topologies $\ic\lambda\tau$ are computed as follows:
\begin{itemize}
\item $\ic 0\tau$ is equal to $\tau$
\item $\ic \lambda\tau$ generated by $\tau$ and all sets of the form
\begin{equation*}
(\alpha, \beta]_\xi^{\tau} := \{x \in X: \alpha < \ell^\xi \rho_\tau x \leq \beta\},
\end{equation*}
for some $-1 \leq \alpha < \beta \leq \Theta$ and some $\xi < \lambda$.
\end{itemize}

\item \label{icardlambdarank}
If $(X,\tau)$ is a scattered space, then  $\rho_{\ic\lambda\tau}=\ell^{\lambda}\circ\rho_\tau$. In particular, the rank function of $\ico \lambda$ is $\ell^\lambda$.
\end{enumerate}
\end{lemma}

\noindent Sets of the form $[\alpha, \beta]_\xi^{\tau}$, $[\alpha, \beta)_\xi^{\tau}$, and $(\alpha, \beta)_\xi^{\tau}$ are defined in the obvious way. In particular, note that $(\alpha, \beta)_0^{\tau} = (\alpha, \beta)^{\tau}$.

\proof
The second claim follows from Lemma \ref{propertiesofdmaps}.\ref{uniquedmap} and Lemma \ref{propertiesofdmaps}.\ref{nextlambdamap}. We use this to prove the first claim by induction. Suppose $\ic \lambda\tau$ is generated by $\tau$ and all sets of the form
\begin{equation*}
(\alpha, \beta]_\xi^{\tau} := \{x \in X: \alpha < \ell^\xi \rho_\tau x \leq \beta\},
\end{equation*}
for $\xi < \lambda$. By definition, $\ic{(\lambda+1)}\tau = \ic 1{(\ic \lambda\tau)}$ is generated by $\ic \lambda\tau$ and all sets of the form
\begin{equation*}
(\alpha, \beta)^{\ic \lambda\tau} := \{x \in X: \alpha < \rho_{\ic \lambda\tau} x < \beta\},
\end{equation*}
but $\rho_{\ic \lambda\tau} = \ell^\lambda \circ \rho_\tau$ by induction hypothesis. So $\ic 1{(\ic \lambda\tau)}$ is generated by all sets of the form
\begin{equation*}
(\alpha, \beta]_\xi^{\tau} := \{x \in X: \alpha < \ell^\xi \rho_\tau x \leq \beta\},
\end{equation*}
for $\xi < \lambda + 1$. The limit case is immediate.
\endproof

The following lemma provides the key relationship between arithmetic and topology for ordinals:
\begin{lemma}[Fern\'andez-Duque \cite{polytopologies}] \label{logdmaps}
Hyperlogarithms 
$$\ell^\xi\colon (\Theta, \ico {\xi+\zeta}) \to (\Theta, \ico\zeta)$$ 
are d-maps.
\end{lemma}

We will make use of the following two lemmata from \cite{AguFer15}:

\begin{lemma}\label{LemmSeparatingNeighborhoods}
Let $(X, \tau)$ be a scattered space and $\lambda$ be an ordinal. Any $x$ in $(X,\ic \lambda\tau)$ has a $\lambda$-neighborhood $U$ such
that whenever $x \neq y \in U$, $\ell^\lambda \rho_0 y < \ell^\lambda \rho_0 x$. 
\end{lemma}

\begin{lemma}\label{nhbase2}
Let $1 < \lambda$ be an additively indecomposable ordinal and $x \in X$ be such that $\rho_{\tau} x = e^\lambda\Theta > 0$. Then for any $\ic\lambda\tau$-neighborhood $V$ of $x$, there exist 
\begin{itemize}
\item a set $U \in \tau$, and 
\item ordinals $\eta < e^\lambda\Theta$ and $\zeta < \lambda$,
\end{itemize} 
such that $V$ contains the set $U \cap (\eta, e^\lambda\Theta]^\tau_{\zeta}$. 
\end{lemma}

For ranks not of the form $e^\lambda\Theta$, we have a more general result, also from \cite{AguFer15}:

\begin{lemma}\label{nhbase3}
Let $(X,\tau)$ be a scattered space. Suppose $0 < \lambda$, $0 < \ell^\lambda\xi$, and $\rho_\tau x = \xi$. Then for any $\ic\lambda\tau$-neighborhood $V$ of $x$, there exist 
\begin{itemize}
\item a set $U \in \tau$, and 
\item a finite partial function $r:\lambda\to\Ord$ such that letting
\[B^X_r(x) = \bigcap_{\zeta \in \text{dom}(r)}(r(\zeta),\ell^\zeta \xi]^\tau_\zeta,\]
\end{itemize} 
we have $U \cap B^X_r(x)\subset V$. 
\end{lemma}	

We conclude this section with a final observation on logarithms:
\begin{lemma}\label{LemmaLeastLog}
Suppose that $\lambda$ is additively indecomposable, $\zeta$ is of the form $e^\lambda \zeta_0$, and $\ell^\lambda \xi < \zeta_0$. Let
\[\eta = \min\{\beta : \ell^\beta\xi<\zeta\}.\]
Then $\eta$ is a successor ordinal or zero.
\end{lemma}
\proof
This is proved by induction on $\xi$.
Suppose towards a contradiction that $\xi$ is least such that $\eta$ is a limit; clearly $\zeta < \xi$. If $\eta$ is additively decomposable, say 
\[\eta = \eta_0 + \omega^\rho > \omega^\rho,\]
then 
\[\ell^\eta\xi = \ell^{\omega^\rho}\ell^{\eta_0}\xi.\]
Now, we must have $\ell^{\eta_0}\xi < \xi$, for otherwise 
\[\xi = \ell^{\eta_0}\xi \leq \ell^{\omega^\rho}\xi = \ell^{\omega^\rho}\ell^{\eta_0}\xi \leq \xi,\]
contradicting the fact $\ell^{\omega^\rho}\ell^{\eta_0}\xi < \zeta$; thus, $\ell^{\eta_0}\xi < \xi$. But then, $\zeta\leq \ell^{\eta_0}\xi$ and
\[\ell^\lambda\ell^{\eta_0}\xi = \ell^{\eta_0 + \lambda}\xi\leq \ell^{\lambda}\xi < \zeta_0,\]
so by the induction hypothesis applied to $\ell^{\eta_0}\xi$, the least $\eta'$ such that 
\[\ell^{\eta'}\ell^{\eta_0}\xi < \zeta\]
is a successor ordinal. However, this ordinal is $\omega^\rho$, which is a contradiction.

Thus $\eta$ is additively indecomposable.
Fern\'andez-Duque and Joosten \cite{hyperations} computed that, letting 
\[\eta^* = \arg\min\{\ell^\nu\xi: 0\leq \nu < \eta\},\]
i.e., letting $\eta^*$ be the least ordinal $\nu$ which minimizes $\ell^\nu \xi$ in $[0,\eta)$, we have
\begin{enumerate}
\item if $0 < \eta^*$, then 
\[\ell^\eta\xi = \ell^\eta\ell^{\eta^*}\xi;\]
\item if $0 = \eta^*$, then
\[\ell^\eta\xi = \sup_{\beta \in [0,\xi)}\big(\ell^\eta\beta+1\big).\]
\end{enumerate}
If $0< \eta^*$, then $\ell^{\eta^*}\xi<\xi$, then one reaches a contradiction as above, using the induction hypothesis on $\ell^{\eta^*}\xi$; thus $0 = \eta^*$. Now, note that $\eta \leq \lambda$, since $\ell^\lambda\xi < \zeta_0\leq\zeta$. In fact, we must have $\eta < \lambda$, by the displayed equation above.
Since $\lambda$ is additively indecomposable, we have
\[\ell^\eta e^\lambda\zeta_0 = e^\lambda\zeta_0 =\zeta <\xi,\]
and thus
\begin{align*}
\ell^\eta \xi 
& = \sup_{\beta \in [0,\xi)}\big(\ell^\eta\beta+1\big)\\
& \geq \ell^\eta\zeta+1\\
&=  \ell^\eta e^\lambda\zeta_0+1\\
& = e^\lambda\zeta_0 + 1,\\
& = \zeta + 1,
\end{align*}
which is again a contradiction. This proves the lemma.
\endproof

\section{$\vec\vartheta$-polytopologies}\label{GenGLP}



In this section, we state our completeness theorem and prove it modulo the \emph{product lemma}, which will be proved in the next section. Let us begin with some motivation by recalling the constructions from \cite{topocompletenessofglp} and \cite{polytopologies}. Let $\mathfrak{X} = (X,\tau)$ be a scattered space; by \cite{AguFer15}, $\gl$ is complete with respect to each topology $\ic\lambda\tau$, with $0<\lambda$, provided $\mathfrak{X}$ is tall enough. Thus, one would attempt to prove completeness of $\glp$ with respect to the polytopology 
\[\{\ic\lambda\tau : 0<\lambda<\Lambda\}.\]
However, this is not a $\glp$-space and thus does not validate the axioms of $\glp$. The idea is then to replace each topology $\ic\lambda\tau$ by a rank-preserving extension and prove completeness for that space. It is not known whether this is possible for arbitrary $\Lambda$. What we will do here is, given a formula $\phi$ consistent with $\glp$, say, with occurrences of modalities $\lambda_0,\hdots, \lambda_n$, and a (tall enough) scattered space $\mathfrak{X} = (X,\tau)$, we produce a sequence of topologies $\tau_0,\hdots, \tau_n$ such that
\begin{enumerate}
\item $(X,\tau_0,\hdots, \tau_n)$ satisfies $\phi$, and
\item each $\tau_i$ is a rank-preserving extension of $\ic{1+\lambda_i}\tau$.
\end{enumerate}

\begin{defi}[$\vartheta$-maximal topology] Let $\vartheta$ be a nonzero ordinal and $(X,\tau)$ be a scattered topological space. We say that $\tau_*$ is a \textit{$\vartheta$-extension} of $\tau$ if 
\begin{enumerate}
\item $\tau \subset \tau_*$, 
\item $\rho_{\tau_*}=\rho_\tau$, and 
\item the identity function $id:(X,\tau)\to(X,\tau_*)$ is continuous at all points $x$ such that $$\ell^\vartheta\rho_\tau(x)=0.$$
\end{enumerate}
We say that $\tau_*$ is an \textit{$\vartheta$-maximal topology} if there are no proper $\vartheta$-extensions of $\tau_*$. 
\end{defi}

In particular, when $\vartheta=1$,  $\vartheta$-maximality coincides with the notion of $\ell$-maximality from \cite{topocompletenessofglp}	. 
If $\vec \vartheta = \{\vartheta_i\colon 0< i < n\}$ is a finite increasing sequence of ordinals, we write 
$$\partial\vec\vartheta := \{\partial\vartheta_{i+1} \colon 0< i < n\},$$
where $\partial\vartheta_{i+1} = - \vartheta_i + \vartheta_{i+1}$ for $0 < i$.  For such a sequence $\vec\vartheta$, we also write $\partial\vartheta_1 = \vartheta_1$.

\begin{defi} Let us call a polytopological space $\mathfrak X = (X, \tau_0, \hdots, \tau_{n})$ a \textit{$\vec \vartheta$-polytopology} over $(X, \tau)$ if $\vec \vartheta = \{\vartheta_i\colon 0 < i \leq n\}$ is increasing and
\begin{enumerate}
\item $\tau_0$ is a $\vartheta_1$-maximal extension of $\tau$,
\item $\tau_{i+1}$ is a $\partial\vartheta_{i+2}$-maximal extension of $\ic{\partial\vartheta_{i+1}}{(\tau_i)}$, for $i +1 < n$, and
\item $\tau_{n} = \ic{\partial\vartheta_{n}}{(\tau_{n-1})}$.
\end{enumerate}
\end{defi}

We remark the following consequence of the definition.
\begin{lemma}
Let $\mathfrak X = (X, \tau_0, \hdots, \tau_{n})$ be a $\vec \vartheta$-polytopology over $(X, \tau)$. Then, for each $i$, 
\[\rho_{\tau_{i}} = \rho_{\ic{\vartheta_i}\tau} = \ell^{\vartheta_i}\circ\rho_\tau.\]
\end{lemma}

Let us refer to the polytopologies considered in \cite{topocompletenessofglp} and \cite{polytopologies} as \emph{$\bg$-polytopologies}. We will not need that notion below, so we do not define them. 
$\vec \vartheta$-polytopologies are weak versions of $\bg$-polytopologies. For example, suppose $\mathfrak X$ is a $\{\omega_1\}$-polytopology over the interval topology on an ordinal. Then $\tau_1$ is a rank-preserving extension of $\ico{\omega_1}$ obtained just by adding sets that would be already included in the corresponding rank-preserving extension of $\ico{\omega_1}$ in any corresponding $\bg$-space of length $\geq \omega_1$ over $\ico 1$. We now state the completeness theorem we shall prove:

\begin{theorem} Let $\vec\vartheta$ be an increasing sequence of nonzero ordinals. Denote by $\glp\upharpoonright \vec\vartheta$ the fragment of $\glp$ whose only modalities appear in $\vec\vartheta$.\label{main3}
\begin{enumerate}
\item (Soundness) All theorems of $\glp\upharpoonright\vec\vartheta$ hold in every $\vec\vartheta$-space.

\item \label{itemcomp} (Completeness) Let $\Lambda \geq \vartheta_{n}$ be a limit ordinal and $(X,\tau)$ be any scattered space of height $\geq e^{\Lambda}1$. Suppose $\varphi$ only contains modalities in $\vec\vartheta$ and is consistent with $\glp\upharpoonright\vec\vartheta$. Then, there is an open subset $U$ of $X$ such that $\varphi$ is satisfied on a $\vec\vartheta$-polytopology over $(U, \ic1\tau)$.
\end{enumerate}
\end{theorem}

The result also holds also for successor ordinals by replacing $e^{\Lambda}1$ with $e^{1 + \Lambda}\omega$. In fact, this general version is what we will prove; the smaller bound in the statement of the theorem follows from the fact that, for limit $\Lambda$,
$$\lim_{\lambda \to \Lambda} e^\lambda\omega 
= \lim_{\lambda \to \Lambda} e^{\lambda+1}1
= \lim_{\lambda \to \Lambda} e^{\lambda}1
\leq \lim_{\lambda \to \Lambda} e^\Lambda 1
= e^\Lambda 1$$
These bounds are sharp (this follows from Lemma \ref{Lemmdmapsiff}). 

Notice that in Theorem \ref{main3}.\ref{itemcomp}, we satisfy the consistent formula on a $\vec\vartheta$-polytopology over a subspace $(X, \ic1\tau)$. We cannot in general replace this with $(X, \ic0\tau)$---consider an ordinal with the initial-segment topology. It is still useful to consider polytopologies of this sort. We will call $\vec\vartheta$-polytopologies over the initial-segment topology \textit{improper}.

In the remainder of this article, we prove Theorem \ref{main3}. Soundness follows from Lemma \ref{nextomegagenbis} below, which in turn follows from Lemma \ref{lmcondOmega}. The proofs are the same as in the case $\vartheta = 1$ from \cite{topocompletenessofglp}. They can also be found in \cite{AguileraThesis}.

\begin{lemma}\label{lmcondOmega}
$(X,\tau)$ is a $\vartheta$-maximal space if, and only if, for all $x \in X$ whose rank $\rho$ is such that $\ell^\vartheta\rho > 0$ and all $V \in \tau$ with $V \subset [0, \rho)_0^\tau$, one of the following holds:
\begin{enumerate}
\item $V \cup \{x\} \in \tau$, or
\item $\rho_\tau(U \cap V) < \rho$ for some $\tau$-neighborhood $U$ of $x$. 
\end{enumerate}
\end{lemma}

\begin{lemma} \label{nextomegagenbis}
Suppose $(X,\tau)$ is $\vartheta$-maximal and $\lambda\geq\vartheta$. Then $\{dA\colon A \subset X\} \subset \ic\lambda\tau$. 
\end{lemma}
It follows from Lemma \ref{nextomegagenbis} that all $\vec\vartheta$-polytopologies are $\glp_n$-spaces, which implies soundness. 

\begin{lemma}[pullback]  \label{Omegapullback}
Suppose $\mathfrak Y = (Y, \mathcal{S}_0, \hdots, \mathcal{S}_{n})$ is a (possibly improper) $\vec\vartheta$-polytopology over $(Y, \sigma)$ and $f\colon (X, \tau) \to (Y, \sigma)$ is a d-map. Then, there exists a $\vec\vartheta$-polytopology $\mathfrak X = (X, \mathcal{T}_0, \hdots, \mathcal{T}_{n})$ over $(X, \tau)$ such that
\begin{equation}\label{eqOmegapullback1}
f\colon (X, \mathcal{T}_i) \to (Y, \mathcal{S}_i)
\end{equation}
is a d-map for each $i \leq n$.
\end{lemma}
\proof
This is essentially the same proof as for the case $\vartheta = 1$ (see \cite[Lemma 8.5]{topocompletenessofglp}). The key point is the following claim:
\begin{claim}
Suppose $\mathcal{T}$ and $\mathcal{S}$ are topologies such that $f\colon (X, \mathcal{T}) \to (Y, \mathcal{S})$ is a d-map and $\mathcal{S}'$ is a $\vartheta$-maximal extension of $\mathcal{S}$, then the topology generated by $\mathcal{T}$ and the family
$$f^{-1}\mathcal{S}' := \{f^{-1}S\colon S \in \mathcal{S}'\}$$
is a $\vartheta$-extension of $\mathcal T$. Moreover, $f\colon (X, \mathcal{T}') \to (Y, \mathcal{S}')$ is a d-map for any $\vartheta$-extension $\mathcal{T}'$ of this topology.
\end{claim}

To see this suffices, suppose the claim holds. Then, letting $\mathcal{T}_0$ be any $\vartheta_1$-maximal extension of the topology given by $\mathcal{T}$ and $f^{-1}\mathcal{S}_0$, we obtain that \eqref{eqOmegapullback1} holds for $i=0$. Inductively, suppose \eqref{eqOmegapullback1} holds for some $i$. By Lemma \ref{propertiesofdmaps}.\ref{nextlambdamap},
$$f\colon(X, \ic{\partial\vartheta_{i+1}}{\mathcal{T}_i}) \to (Y, \ic{\partial\vartheta_{i+1}}{\mathcal{S}_i})$$
is a d-map. If $i+1 = n$, then we are done; otherwise, by definition, $\mathcal{S}_{i+1}$ is a $\partial\vartheta_{i+2}$-maximal extension of $\ic{\partial\vartheta_{i+1}}{\mathcal{S}_i}$, whereby the claim yields that \eqref{eqOmegapullback1} holds for $i+1$ if we set $\mathcal{T}_{i+1}$ to be some $\partial\vartheta_{i+2}$-maximal extension of the topology given by $\ic{\partial\vartheta_{i+1}}{\mathcal{T}_i}$ and $f^{-1}\mathcal{S}_{i+1}$. Hence, it suffices to prove the claim.

\proof[Proof of the claim]
Let $\mathcal R$ be the topology given by $\mathcal{T}$ and $f^{-1}\mathcal{S}'$. Using the fact that
\begin{equation} \label{eqOmegapullback2}
f\colon (X, \mathcal{T}) \to (Y, \mathcal{S}) \text{ is a d-map},
\end{equation}
it is not hard to see that $f\colon (X, \mathcal R) \to (Y, \mathcal S')$ is also a d-map. By definition, $\mathcal T \subset \mathcal R$, whence $\mathcal R$ is a rank-preserving extension of $\mathcal T$. Let $x \in X$ be such that $\ell^{\vartheta} \rho_{\mathcal{T}}x = 0$. We need to show that $id\colon (X, \mathcal T) \to (X, \mathcal R)$ is continuous at $x$. This follows from the fact that $\mathcal{S}'$ is a $\vartheta$-extension of $\mathcal S$: for any $\mathcal R$-neighborhood of $x$ of the form $f^{-1}V$, we have that $V$ is a $\mathcal S'$-neighborhood of $f(x)$ and $\ell^{\vartheta} \rho_{\mathcal{S}'}f(x) = 0$, so that $f(x) \in U$ for some $U \in \mathcal S$ with $U \subset V$. Therefore, $x \in f^{-1}U \subset f^{-1}V$ and $f^{-1}U \in \mathcal T$, by \eqref{eqOmegapullback2}.

Now let $\mathcal T'$ be any $\vartheta$-extension of $\mathcal R$. Clearly, $f\colon (X, \mathcal{T}') \to (Y, \mathcal{S}')$ is continuous and pointwise discrete. Suppose towards a contradiction that $x \in X$ and $W \in \mathcal T'$ witness a failure of $f$ being open. Let 
$$\rho := \rho_{\mathcal T}x = \rho_{\mathcal T'}x = \rho_{\mathcal S'}f(x).$$
Note that we must have $0 < \ell^\vartheta \rho$. Without loss of generality, we may assume $\rho$ is the least possible rank of a counterexample and $W$ contains no other point of rank $\geq \rho$, so that $W = W_0 \cup \{x\}$, for some $W_0 \in \mathcal{T}$ with $W_0\subset [0, \rho)^{\mathcal T}_0$. We will arrive at a contradiction using Lemma \ref{lmcondOmega}: since $f$ is rank-preserving, we have that $\ell^\vartheta \rho > 0$, $f(W_0) \in \mathcal{S}'$, and $f(W_0)\subset [0, \rho)^{\mathcal S'}_0$. Hence, by Lemma \ref{lmcondOmega}, one of the following holds:
\begin{enumerate}
\item $f(W_0) \cup \{f(x)\} \in \mathcal{S}'$, or
\item \label{pbO2} $\rho_{\mathcal{S}'}(U \cap f(W_0)) < \rho$ for some $\mathcal S'$-neighborhood $U$ of $f(x)$. 
\end{enumerate}
It must be the second one that holds, for $f(W_0) \cup \{f(x)\} = f(W)$ is not $\mathcal{S}'$-open by hypothesis. 
Observe that $f^{-1}U \cap W$ is a $\mathcal T'$-neighborhood of $x$ and thus contains points with rank of every ordinal up to, and including, $\rho$. Because $W$ contains only one point of rank $\rho$ and $f$ is rank-preserving,
$$(f^{-1}U \cap W) \setminus f^{-1}(x) = f^{-1}U \cap W_0.$$
It follows that the set
\[f^{-1}U \cap W_0\]
contains points with rank of every ordinal up to, but not including, $\rho$. However, this is impossible by \ref{pbO2} above, because $\rho$ is a limit ordinal. This finishes the proof of the claim and the lemma.
\endproof

Lemma \ref{Omegapullback} is still true in the degenerate case $\sigma = \ico 0$. In this case, notice that $\ico 0$ is already $\vartheta$-maximal for every $\vartheta$, for there is only one point of each rank. The proof of Theorem \ref{LemmTopoO} below is postponed to the next section.
\begin{theorem}[Product Lemma]\label{LemmTopoO} 
Assume $\varsigma$ is a nonzero additively indecomposable ordinal, $([1, \kappa], \mathcal{T}_0, \hdots, \mathcal{T}_{n})$ is a $\vec\vartheta$-polytopology over $\ico\varsigma$, and $([1, \lambda], \mathcal{S}_0, \hdots, \mathcal{S}_{n})$ is a $\vec\vartheta$-polytopology over $\ico 0$. Suppose moreover that
\[\kappa < e^{\varsigma+\vartheta_n}\omega\]
and that
\[\lambda < e^{\vartheta_n}\omega.\]
Fix $\{\kappa_0, \hdots, \kappa_m\}$ a finite subset of $[1, \kappa]$. Let \begin{equation*}
\xi = \text{ the unique ordinal equal to } \ell^\varsigma[1, \kappa_m].
\end{equation*}
Let $\Theta = e^\varsigma(\xi + (-1 + \lambda))< e^{\varsigma+\vartheta_n}\omega$ and define $X_\uparrow = [1,\Theta] \cap [\xi,\infty)_\varsigma$ and $X_\downarrow = [1,\Theta]\cap[0,\xi)_\varsigma$. Then, there exist:

\begin{enumerate}

\item \label{LemmTopoO2} A $\vec\vartheta$-polytopology $([1, \Theta], \mathcal R_0, \hdots, \mathcal R_n)$ over $\ico \varsigma$.
 
\item \label{LemmTopoO3}  Functions $\pi_0 \colon [1, \Theta] \to [1, \kappa]$ and $\pi_1 \colon [1, \Theta] \to [1, \lambda]$ such that:
\begin{itemize}

\item $\pi_0 \colon (X_\downarrow, \mathcal R_i)$ $\to ([1,\kappa], \mathcal T_i)$ is a surjective d-map for each $i$;

\item $\pi_1 \colon (X_\uparrow, \mathcal R_i)$ $\to ([1,\lambda], \mathcal S_i)$ is a surjective d-map for each $i$;

\item $X_\uparrow \subset d_{\mathcal{R}_0}\pi_0^{-1}(\kappa_i)$ for any $i<m$;

\item $\pi_1 =1 + (-\xi + \ell^\varsigma)$;

\item the polytopology $(\mathcal{R}_0,\hdots, \mathcal{R}_n)$, when restricted to $X_\uparrow$, is the one obtained from Lemma \ref{Omegapullback} by pulling back via $\pi_1$;

\item $\pi_1^{-1}(\{\lambda\}) = \{\Theta\}$.
\end{itemize}
\end{enumerate}
\end{theorem}

Theorem \ref{LemmTopoO} is the main new ingredient of our proof. With it, we can adapt the usual proofs to obtain completeness. First, we need an embedding lemma:

\begin{lemma}\label{tm2bis}
Let $(T, <_0, <_1, ..., <_n)$ be a finite J$_n$-tree with root $r$ and $\vec \vartheta$ be an increasing $n$-sequence of nonzero ordinals. Then, for any $\varsigma > 0$, there exist 
\begin{itemize}

\item a $\vec\vartheta$-polytopology $([1, \Theta], \mathcal{T}_0, \hdots, \mathcal{T}_{n})$ over $([1, \Theta], \ico \varsigma)$ such that $\Theta < e^{\varsigma + \vartheta_n}\omega$; and

\item a surjective J$_n$-map $f: ([1, \Theta], \mathcal{T}_0, \hdots, \mathcal{T}_{n})\to T$ such that $f^{-1}(r)=\{\Theta\}$.
\end{itemize}
\end{lemma}

\proof
The proof is by induction on $n$. The base case follows from \cite[Theorem 6.11]{AguFer15}, so we assume that the result holds for all $m < n$ and proceed by a subsidiary double induction on
\begin{enumerate}
\item $\varsigma$, which we decompose as $\varsigma_0 + \omega^\rho$, and
\item $hgt_0(T)$, the height of $<_0$,
\end{enumerate}
in that order.
Let $hgt_i(T)$ be the height of $<_i$. We need to consider various cases:

\paragraph{Case I:} $\omega^\rho < \varsigma$. By the induction hypothesis (for $\varsigma$) applied to $\omega^\rho$, there are:
\begin{itemize}
\item a $\vec\vartheta$-polytopology $([1, \Theta_0], \mathcal{S}_0, \hdots, \mathcal{S}_{n})$ over $([1, \Theta_0], \ico {\omega^\rho})$ such that $\Theta_0 
< e^{\omega^\rho + \vartheta_n}\omega$; and

\item a surjective J$_n$-map $g: ([1, \Theta_0], \mathcal{S}_0, \hdots, \mathcal{S}_{n})\to T$ such that $g^{-1}(r)=\{\Theta_0\}$.
\end{itemize}

\noindent By Lemma \ref{logdmaps}, 
$$\ell^{\varsigma_0}\colon ([1, e^{\varsigma_0}\Theta_0], \ico{\varsigma}) \to ([1, \Theta_0], \ico {\omega^\rho})$$
is a d-map, whence by Lemma \ref{Omegapullback}, there is a $\vec\vartheta$-polytopology 
$$\mathfrak X = ([1, e^{\varsigma_0}\Theta_0], \mathcal{T}_0, \hdots, \mathcal{T}_n)$$
over $([1, e^{\varsigma_0}\Theta_0], \ico{\varsigma})$ such that 
$$\ell^{\varsigma_0}\colon ([1, e^{\varsigma_0}\Theta_0], \mathcal{T}_i) \to ([1, \Theta_0], {\mathcal{S}_i})$$
is a d-map for each $i \leq n$. It follows that $f := g \circ \ell^{\varsigma_0}$ is a surjective J$_n$-map from $\mathfrak X$ onto $T$. Finally, 
$$\Theta:= e^{\varsigma_0}\Theta_0 
< e^{\varsigma_0 + \omega^\rho + \vartheta_n}\omega 
= e^{\varsigma + \vartheta_n}\omega.$$

\paragraph{Case II:} $\varsigma$ is additively indecomposable and $0 = hgt_0(T)$, so that $<_0 = \varnothing$. Let 
$$\vec\vartheta^* = -\vartheta_1 + \vec\vartheta\upharpoonright[2, n] = (-\vartheta_1 + \vartheta_2,-\vartheta_1 + \vartheta_3,\hdots, -\vartheta_1 + \vartheta_n).$$ By induction hypothesis (applied to $n$), there are:
\begin{itemize}
\item a $\vec\vartheta^*$-polytopology $([1, \Theta_0], \mathcal{S}_1, \hdots, \mathcal{S}_{n})$ over $([1, \Theta_0], \ico {\vartheta_1})$ such that 
$$\Theta_0 < e^{\vartheta_1 + (-\vartheta_1 + \vartheta_n)}\omega = e^{\vartheta_n}\omega;$$

\item a surjective J$_{n-1}$-map $g: ([1, \Theta_0], \mathcal{S}_1, \hdots, \mathcal{S}_{n})\to (T, <_1,\hdots, <_n)$ such that $g^{-1}(r)=\{\Theta_0\}$.
\end{itemize}
Note that each $\mathcal{S}_i$ is a rank-preserving extension of $\ico{\vartheta_i}$.
Thus, a simple computation shows that $([1, \Theta_0], \ico0, \mathcal{S}_1, \hdots, \mathcal{S}_{n})$ is a(n improper) $\vec\vartheta$-polytopology. Clearly,
$$\ell^{\varsigma}\colon ([1, e^{\varsigma}\Theta_0], \ico\varsigma) \to ([1, \Theta_0], \ico 0)$$
is a d-map. By Lemma \ref{Omegapullback}, there exists a $\vec\vartheta$-polytopology $\mathfrak X$ of the form $([1, e^{\varsigma}\Theta_0], \mathcal{T}_0, \hdots, \mathcal{T}_{n})$ over $([1, e^{\varsigma}\Theta_0], \ico \varsigma)$ such that
\begin{equation*}
\ell^{\varsigma}\colon  ([1, e^{\varsigma}\Theta_0], \mathcal{T}_i) \to ([1, \Theta_0], \mathcal{S}_i)
\end{equation*}
is a d-map for each $i \leq n$. Let $\Theta := e^{\varsigma}\Theta_0 < e^{\varsigma + \vartheta_n}\omega$. Let $f:= g \circ \ell^{\varsigma}$. We have that
\begin{equation} \label{eqjmapemb}
f\colon  ([1, \Theta], \mathcal{T}_1, \hdots \mathcal{T}_n) \to (T, <_1, \hdots, <_n) \text{ is a J$_{n-1}$-map}.
\end{equation}

We claim \eqref{eqjmapemb} holds for the full space $\mathfrak X$ and the structure $(T, <_0, \hdots, <_n)$, i.e., $f$ is already a J$_n$-map: condition ($j_1$) is given by definition; ($j_2$) is satisfied trivially since the topology induced by $<_0$ is discrete; ($j_3$) and ($j_4$) hold for $k\neq 0$ because of \eqref{eqjmapemb} and for $k = 0$ because $r$ is the sole hereditary $1$-root of $T$. 

\paragraph{Case III:} $\varsigma$ is additively indecomposable and $0 < hgt_0(T)$. 
Let $r_1, \hdots, r_m$ be all $<_0$-successors of $r$ that are hereditary $1$-roots and (following earlier notation) let ${\ll_0}(r_i)$ denote the generated subtrees. Also let ${\ll_1}(r)$ denote the subtree consisting of all nodes that are $<_0$-incomparable with $r$ (i.e., the $<_0$-roots). By induction hypothesis (applied to $hgt_0(T)$), there exist $\vec\vartheta$-polytopologies
$$\mathfrak{X}_i = ([1, \kappa_i], \mathcal{T}_0^i,\hdots, \mathcal{T}_n^i)$$
over $([1, \kappa_i], \ico\varsigma)$ and surjective J$_n$-maps 
$$f_i\colon\mathfrak X_i \to ({\ll_0}(r_i), <_0,\hdots, <_n)$$ 
for $0<i$. Let 
\[\kappa = \kappa_1 + \hdots + \kappa_m\] 
and 
$$\mathfrak X = ([1, \kappa], \mathcal{T}_0,\hdots\mathcal{T}_n)$$ 
be the topological sum. 
We also denote by 
$$f^*: \mathfrak{X} \to \left(\bigcup_{0<i\leq n}{\ll_0}(r_i), <_0, \hdots, <_n\right)$$ 
the sum of the functions $f_i$. 
We may define an improper $\vec\vartheta$-polytopology 
$$\mathfrak Y' = ([1,\lambda], \ico 0, \mathcal{S}_1,\hdots, \mathcal{S}_n)$$ 
over $\ico 0$ as in Case II in such a way that there is a J$_{n-1}$-map 
$$g\colon ([1,\lambda],\mathcal{S}_1,\hdots,\mathcal{S}_n) \to ({\ll_1}(r), <_1, \hdots, <_n)$$ 
such that, letting $f_0 = \ell^\varsigma\circ g$ and 
$$\mathfrak Y = [1, e^\varsigma\lambda]$$
then
$$f_0\colon \mathfrak Y \to ({\ll_1}(r), <_0, ..., <_n)$$ 
is a J$_n$-map. In fact, if we instead define
$$f_0(\zeta) = g\Big(1+ \big(-\xi + \ell^\varsigma(\zeta)\big)\Big)$$ then
$$f_0: \Big(\big[1,e^\varsigma(\xi + (-1 + \lambda))\big] \cap \big[\xi,\infty\big)_\varsigma, \hat{\mathcal{S}_0}, \hdots, \hat{\mathcal{S}}_n\Big) \to \big({\ll_1}(r),<_0, \hdots, <_n\big)$$ 
is also a J$_n$-map, where $\hat{\mathcal{S}_0}, \hdots, \hat{\mathcal{S}}_n$ is the polytopology obtained from $\ico 0, \mathcal{S}_1,\hdots, \mathcal{S}_n$ via Lemma \ref{Omegapullback}. Let
\[\Theta = e^\varsigma(\xi + (-1 + \lambda))\]
and write $X_\uparrow = [1,\Theta] \cap [\xi,\infty)_\varsigma$ and $X_\downarrow = [1,\Theta]\cap[0,\xi)_\varsigma$. 
By the Product Lemma, there are:
\begin{enumerate}

\item  A $\vec\vartheta$-polytopology $([1, \Theta], \mathcal R_0, \hdots, \mathcal R_n)$ over $\ico \varsigma$.
 
\item  Functions $\pi_0 \colon [1, \Theta] \to [1, \kappa]$ and $\pi_1 \colon [1, \Theta] \to [1, \lambda]$ such that:
\begin{itemize}

\item $\pi_0 \colon (X_\downarrow, \mathcal R_i)$ $\to ([1,\kappa], \mathcal T_i)$ is a surjective d-map for each $i$;

\item $\pi_1 \colon (X_\uparrow, \mathcal R_i)$ $\to ([1,\lambda], \mathcal S_i)$ is a surjective d-map for each $i > 0$;

\item $X_\uparrow \subset d_{\mathcal{R}_0}\pi_0^{-1}(\kappa_i)$ for any $i<m$;

\item $\pi_1 =1 + (-\xi + \ell^\varsigma)$;

\item the polytopology $(\mathcal{R}_0,\hdots, \mathcal{R}_n)$, when restricted to $X_\uparrow$, is the one obtained from Lemma \ref{Omegapullback} by pulling back via $\pi_1$;

\item $\pi_1^{-1}(\{\lambda\}) = \{\Theta\}$.
\end{itemize}
\end{enumerate}
Thus, when restricted to $X_\uparrow$, the topologies $\mathcal{R}_i$ and $\hat{\mathcal{S}_i}$ coincide, as do the mappings $\pi_1$ and $f_0$, and so
$$f_0: (X_\uparrow, \mathcal R_0, \hdots, \mathcal R_n) \to ({\ll}(r),<_0,\hdots, <_n)$$
is a J$_n$-map.
We define a function
$$f: ([1, \Theta], \mathcal R_0, \hdots, \mathcal R_n) \to (T,<_0,\hdots, <_n)$$
given by:

\begin{equation*}
f(x) =
\left\{
	\begin{array}{ll}
		f^*(\pi_0(x))  & \mbox{if } x \in X_\downarrow\\
		f_0  & \mbox{if } x \in X_\uparrow 
	\end{array}
\right.
\end{equation*}
Since $X_\uparrow$ and $X_\downarrow$ are $\ico{\varsigma+1}$-clopen and $1\leq\partial\vartheta_1$, it follows that $X_\uparrow$ and $X_\downarrow$ are $\mathcal{R}_i$-clopen for all $0<i$.
The facts that $f_0$ and $f^*$ are J$_n$-maps and that the projection $\pi_0$ is a d-map yield condition ($j_1$), as well as conditions ($j_2$)--($j_4$) for $0<i$. We verify the remaining ones:

\begin{enumerate}
\item[($j_2$)] Let $U$ be $\mathcal{R}_{0}$-open. If $U \subset X_\downarrow$, then $f(U)$ is $<_0$-open, as $\pi_0\circ f^*$ is a J$_n$-map. If $U \cap X_\uparrow \neq \varnothing$, then we claim $f(U)$ is $<_0$-open in $T$. Indeed, since $X_\uparrow \subset d_{\mathcal{R}_0}\pi_0^{-1}(\kappa_i)$ for any $0<i\leq m$, then there are ordinals $\xi_0, \hdots, \xi_m \in U$ such that $\pi_0(\xi_i) = \kappa_i$ for each $i$. But then $\pi_0(U\cap X_\downarrow)$ contains a neighborhood $U_i$ of each $\kappa_i$ and by choice of $f^*$, $f^*(U_i) = {\ll_0}(r_i)$.

\item[($j_3$)] Any hereditary $1$-root $x$ is either $r$ or in some $T_i$. In the former case, $f^{-1}({\ll_0}(r))$ and $f^{-1}(\{r\}\cup {\ll_0}(r)) \in \tau_0$ equal $[1, \Theta)$ and $[1, \Theta]$, respectively. In the latter case, the result follows from the continuity of $\pi_0$ and the fact that $f^*$ is a J$_n$-map.

\item[($j_4$)] Again, $f^{-1}(\{r\}) = \{\Theta\}$ and for any hereditary $1$-root $x \neq r$, $f^{-1}(x) = \pi_0^{-1}f^{*-1}(x)$ is discrete because $f^{*-1}(x)$ is discrete and $\pi_0$ is pointwise discrete.
\end{enumerate}

Therefore, $f$ is a indeed a surjective J$_n$-map.\\

\noindent Since we have considered all cases, the lemma follows.
\endproof

We can now finish the proof of Theorem \ref{main3}. Let $\Lambda \geq \vartheta_{n}$ and $(X,\tau)$ be any scattered space of height $\geq e^{1+\Lambda}\omega$. Suppose $\varphi$ only contains modalities in $\vec\vartheta$ and is consistent with $\glp\upharpoonright\vec\vartheta$. We need to show that $\varphi$ is satisfied on a $\vec\vartheta$-polytopology over $(X, \ic1\tau)$.

We use Lemmata \ref{to comp 4} and  \ref{tm2bis} to find 
\begin{enumerate}
\item the canonical tree $T$ for $\varphi$ with root $r$, 
\item a $\vec\vartheta$-polytopology $([1,\Theta], \mathcal{T}_0,\hdots, \mathcal{T}_n)$ over $([1,\Theta], \ico 1)$ such that $\Theta< e^{1 + \vartheta_n}\omega$, and
\item a surjective J$_n$-map $f:([1,\Theta], \mathcal{T}_0,\hdots, \mathcal{T}_n) \to T$ such that $f^{-1}(r) = \{\Theta\}$.
\end{enumerate}
By Lemma \ref{to comp 4},
\[([1,\Theta], \mathcal{T}_0,\hdots, \mathcal{T}_n)\not\models\varphi.\]
Now, by assumption, $(X, \tau)$ is a scattered space such that 
$$e^{1+\vartheta_n}\omega  \leq e^{1+\Lambda}\omega \leq hgt(X,\tau).$$ 
Choose $k<\omega$ such that
$$\Theta< e^{1 + \vartheta_n}k$$
and choose a point $x \in X$ of rank $e^{1 + \vartheta_n}k$. Let $U$ be a neighborhood of $x$. Without loss of generality, $(U,\tau)$ is a scattered space of height $e^{1 + \vartheta_n}k + 1$. 
By Lemma \ref{Omegapullback}, there is a $\vec\vartheta$-polytopolgy 
\[\mathfrak{U} = (U,\mathcal{S}_0,\hdots, \mathcal{S}_n)\]
over $(U, \ic 1 \tau)$ such that
\[\rho_\tau: (U,\mathcal{S}_i)\to ([1,e^{1 + \vartheta_n}k], \mathcal{T}_i)\]
is a d-map for each $i\leq n$.
By Lemmata \ref{to comp 4} and \ref{jmapcomposition},  $\varphi$ is then satisfiable in $\mathfrak U$. 
This completes the proof of Theorem \ref{main3}.


\section{Proof of the Product Lemma}\label{SectLifting}
For convenience, we restate the lemma:
\begin{theorem}[Product Lemma] 
Assume $\varsigma$ is a nonzero additively indecomposable ordinal, $([1, \kappa], \mathcal{T}_0, \hdots, \mathcal{T}_{n})$ is a $\vec\vartheta$-polytopology over $\ico\varsigma$, and $([1, \lambda], \mathcal{S}_0, \hdots, \mathcal{S}_{n})$ is a $\vec\vartheta$-polytopology over $\ico 0$. Suppose moreover that
\[\kappa < e^{\varsigma+\vartheta_n}\omega\]
and that
\[\lambda < e^{\vartheta_n}\omega.\]
Fix $\{\kappa_0, \hdots, \kappa_m\}$ a finite subset of $[1, \kappa]$. Let \begin{equation*}
\xi = \text{ the unique ordinal equal to } \ell^\varsigma[1, \kappa_m].
\end{equation*}
Let $\Theta = e^\varsigma(\xi + (-1 + \lambda))< e^{\varsigma+\vartheta_n}\omega$ and define $X_\uparrow = [1,\Theta] \cap [\xi,\infty)_\varsigma$ and $X_\downarrow = [1,\Theta]\cap[0,\xi)_\varsigma$. Then, there exist:

\begin{enumerate}

\item  A $\vec\vartheta$-polytopology $([1, \Theta], \mathcal R_0, \hdots, \mathcal R_n)$ over $\ico \varsigma$.
 
\item  Functions $\pi_0 \colon [1, \Theta] \to [1, \kappa]$ and $\pi_1 \colon [1, \Theta] \to [1, \lambda]$ such that:
\begin{itemize}

\item $\pi_0 \colon (X_\downarrow, \mathcal R_i)$ $\to ([1,\kappa], \mathcal T_i)$ is a surjective d-map for each $i$;

\item $\pi_1 \colon (X_\uparrow, \mathcal R_i)$ $\to ([1,\lambda], \mathcal S_i)$ is a surjective d-map for each $i$;

\item $X_\uparrow \subset d_{\mathcal{R}_0}\pi_0^{-1}(\kappa_i)$ for any $i<m$;

\item $\pi_1 =1 + (-\xi + \ell^\varsigma)$;

\item the polytopology $(\mathcal{R}_0,\hdots, \mathcal{R}_n)$, when restricted to $X_\uparrow$, is the one obtained from Lemma \ref{Omegapullback} by pulling back via $\pi_1$;

\item $\pi_1^{-1}(\{\lambda\}) = \{\Theta\}$.
\end{itemize}
\end{enumerate}
\end{theorem}

A slight change in notation will make the proof easier: instead of starting with an ordinal $\kappa$ and a finite subset $\{\kappa_1, \hdots, \kappa_m\}$ of $[1, \kappa]$, we will start with a set of ordinals $\{\kappa_1, \hdots, \kappa_m\}$ and define $\kappa$ as their sum. 

So let $\varsigma$ be an additively indecomposable ordinal (possibly equal to $1$). Let $\lambda, \kappa_1, \hdots, \kappa_m$ be nonzero ordinals such that $\kappa_1 \leq \dots \leq \kappa_m$ and write 
$$\kappa = \kappa_1 + \hdots + \kappa_m.$$ 
We will often speak of ordinals $\mod m$ and of ``$\kappa_{\iota \mod m}$.'' In particular, 
\[\kappa_0 = \kappa_{m \mod m}\]
is the greatest among all of $\kappa_1$, $\kappa_2, \hdots, \kappa_m$.
By convention, we assume that $\omega \equiv 0 \mod m$.

Fix
  $\vec\vartheta$-polytopologies 
  $([1, \kappa], \mathcal{T}_0, \hdots, \mathcal{T}_{n})$ and $([1, \lambda], \mathcal{S}_0, \hdots, \mathcal{S}_{n})$
  over $\ico\varsigma$ and $\ico0$, respectively and
\begin{equation*}
\xi = \text{ the unique ordinal equal to } \ell^\varsigma[1, \kappa_m].
\end{equation*}
Note that $\xi$ exists and is a successor ordinal. Write $\zeta$ for the predecessor of $\xi$. We have $\kappa_m < e^\varsigma\xi$, for otherwise $e^\varsigma\xi$ belongs to the interval $[1, \kappa_m]$ and so $\xi \in \ell^\varsigma[1,\kappa_m]$, contrary to its definition.
Now, condition \eqref{FixPoints} states that
\[e^\varsigma(\zeta + 1) = \lim_{\alpha\to\varsigma}e^\alpha(e^\varsigma\zeta + 1).\]
Since $e^\varsigma\zeta\leq\kappa_m < e^\varsigma\xi$, it follows that 
$$e^\alpha(e^\varsigma\zeta + 1) \leq \kappa_m$$ 
for some greatest ordinal $\alpha$, which we will denote by $\nu$, i.e., 
\begin{equation*}
e^\nu(e^\varsigma\zeta + 1) \leq \kappa_m < e^{\nu + 1}(e^\varsigma\zeta + 1).
\end{equation*}

The proof of the Product Lemma is distributed among a series of lemmata and definitions throughout this section.

\begin{defi} \label{DefCS}
We define the \textit{characteristic sequence for $\varsigma$}, $\csp{\varsigma}$ as follows:
\begin{itemize}
\item If $\varsigma = 1$, then $\csp{\varsigma}$ is the sequence $\{\varsigma_\iota\}_{\iota < \omega}$ with constant value $0$.

\item If $\varsigma$ is a limit ordinal, then $\csp{\varsigma}$ is the sequence $\{\varsigma_\iota\}_{\iota < \varsigma}$ given by $\varsigma_\iota = \max\{\iota, \nu\cdot\iota\}$.
\end{itemize}
\end{defi}

\begin{remark}
The case $\varsigma = 1$ needs to be considered separately. It corresponds to the notion of d-products from \cite{topocompletenessofglp}. In our context, the Product Lemma will be sufficient, but it is also possible to ignore the case $\varsigma = 1$ of the Product Lemma and instead employ the d-product construction directly (cf. \cite{AguileraThesis}).
\end{remark}

Clearly, if $\varsigma$ is multiplicatively indecomposable, then the sequence $\cs$ is cofinal in $\varsigma$.
We shall assume for notational simplicity that $\varsigma$ is multiplicatively indecomposable, as opposed to simply additively indecomposable, hereafter. If it were not, however, multiplicatively indecomposable, then it must still be a limit of points in the sequence $\cs$. 
If $\varsigma$ were not multiplicatively indecomposable, then the following proof goes through if one replaces the sequence $\cs$ with the appropriate initial segment throughout.

\subsection{The partition}\label{subsectionGLPPartition}
For the rest of the section, we shall assume that $\nu$ is nonzero; this will avoid various awkward case distinctions. If $\nu$ were in fact $0$, or even finite, then the computations carried out below would in fact be simpler, and similar to the situation in \cite{AguFer15}.
\begin{defi}[$\Theta, X_\uparrow, X_\downarrow$]\label{cellpartition}
For each $\iota < \varsigma$, we take $\kappa_\iota$ to mean the unique $\kappa_i$ such that $\iota \equiv i \mod m$. In particular, $\kappa_0 = \kappa_m \geq \kappa_i$ for any $i$. Set:

\begin{itemize}
\item $\alpha_0 = 1$;
\item $\alpha_{\iota + 1} = \beta_\iota + 1$;
\item $\alpha_{\iota} = e^{\varsigma_\iota}(e^\varsigma\zeta + 1)$, at limit stages; and
\item $\beta_\iota = 1 + e^{\varsigma_\iota}(\kappa_0 + \kappa_\iota)$
\end{itemize}

Let $X_\iota = [\alpha_\iota, \beta_\iota]$, $Y_\iota = X_\iota \cap [0, \kappa_0]_{\varsigma_\iota}$, and $Z_\iota = X_\iota \setminus Y_\iota$. 
We will soon prove (Lemma \ref{lift2}) that the family $\{X_\iota\colon\iota < \varsigma\}$ partitions $e^\varsigma\xi$:
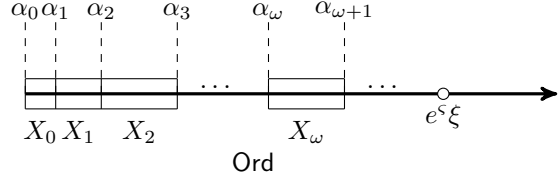
\begin{center}
\begin{tikzpicture}[
    axis/.style={very thick, ->, >=stealth'},
    ll/.style={thick},
    dline/.style={dashed, thin}, 
    pile/.style={thick, ->, >=stealth', shorten <=2pt, shorten
    >=2pt},
    wdot/.style={ 
    circle,
    fill=white,
    draw,
    outer sep=0pt,
    inner sep=1.5pt
  }
    ]
\draw[axis] (0,0) -- (7,0);
\node at (3,-0.9) {$\Ord$};

\foreach \x in {0, 0.4, 1, 2}
{
	\draw[dline] (\x,0) -- (\x,1);
	\draw        (\x,-0.2) -- (\x,0.2);
}

	\node at (0,1.1) {$\alpha_0$};
	\node at (0.4,1.1) {$\alpha_1$};
	\node at (1,1.1) {$\alpha_2$};
	\node at (2,1.1) {$\alpha_3$};	
	\node at (2.5,0.1) {$\hdots$};
	
\foreach \x in {3.2, 4.2}
{
	\draw[dline] (\x,0) -- (\x,1);
	\draw        (\x,-0.2) -- (\x,0.2);
}
	\node at (3.2,1.1) {$\alpha_\omega$};
	\node at (4.2,1.1) {$\alpha_{\omega+1}$};	
	\node at (4.7,0.1) {$\hdots$};
	
	\foreach \x in {-0.2, 0.2}
{
	\draw        (0,\x) -- (2,\x);
	\draw        (3.2,\x) -- (4.2,\x);
}
	\node at (0.2,-0.5) {$X_0$};
	\node at (0.7,-0.5) {$X_1$};
	\node at (1.5,-0.5) {$X_2$};
	\node at (3.7,-0.5) {$X_\omega$};
	
	\node[wdot] at (5.5,0) {};
	\node at (5.5,-0.3) {$e^\varsigma\xi$};
	
\end{tikzpicture}
\end{center}
\end{defi}

\begin{defi}
As in the statement of the theorem, we define:
\begin{align*}
\lift 
&= \text{ least } \alpha \text{ such that } 1 + (-\xi + \ell^\varsigma\alpha) = \lambda\\
&= e^\varsigma(\xi + (-1 + \lambda))\\
&= e^\varsigma( \ell^\varsigma[1,\kappa] + (-1 + \lambda)).
\end{align*}
Observe that if $\lambda < e^{\vartheta_n}\omega$ and $\kappa < e^{\varsigma+\vartheta_n}\omega$, then
\[e^\varsigma( \ell^\varsigma[1,\kappa] + (-1 + \lambda)) < e^{\varsigma+\vartheta_n}\omega,\]
as desired.
We also set:
\begin{enumerate}
\item $X_\downarrow = [1, \lift] \cap [0, \xi)_\varsigma$; \item $X_\uparrow = [1, \lift] \cap [\xi, \infty)_\varsigma$. 
\end{enumerate} 
See the following picture:
\end{defi}

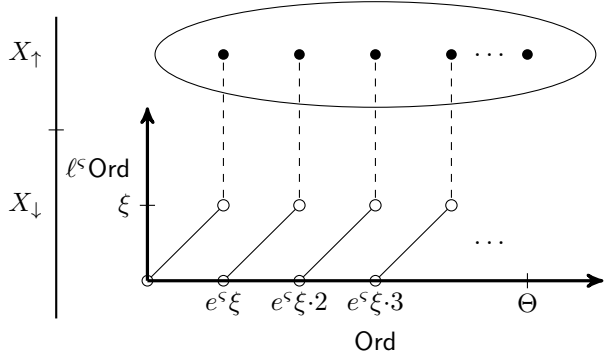
\begin{center}
\begin{tikzpicture}[
    axis/.style={very thick, ->, >=stealth'},
    ll/.style={thick},
    dline/.style={dashed, thin}, 
    pile/.style={thick, ->, >=stealth', shorten <=2pt, shorten
    >=2pt},
    wdot/.style={ 
    circle,
    fill=white,
    draw,
    outer sep=0pt,
    inner sep=1.5pt
  }
    ]
\node at (-1.6,1) {$X_\downarrow$};    
\node at (-1.6,3) {$X_\uparrow$};    
\draw (-1.3,2) -- (-1.1,2); 
\draw[ll] (-1.2,-0.5) -- (-1.2,3.5);
\draw (3,3) ellipse (2.9 and 0.7);
    
\draw[axis] (0,0) -- (0,2.3); 
\node at (-0.65,1.5) {$\ell^{\varsigma}\Ord$};
\draw[axis] (0,0) -- (6,0); 
\node at (3,-0.8) {$\Ord$};
    
\node at (-0.3,1) {$\xi$}; 
\draw (-0.1,1) -- (0.1,1); 

\foreach \x in {0, 1, 2, 3}
{
	\draw (\x,0) circle (2pt);
	\draw (\x+0.05,0.05) -- (\x+1,1);
	\node[wdot] at (\x+1,1) {};
	\draw[dline] (\x+1,1.1) -- (\x+1,2.9);
	\fill (\x+1,3) circle [radius=2pt];
}
	\node at (4.5,0.5) {$\hdots$};
	\node at (4.5,3) {$\hdots$};
	
	\fill (5,3) circle [radius=2pt];
	\draw (5,-0.1) -- (5,0.1);
	\node at (5,-0.3) {$\Theta$};
	
	\node at (1,-0.3) {$e^\varsigma\xi$};
\foreach \x in {2,3}
{
	\node at (\x,-0.3) {$e^\varsigma\xi\cdot$\x};
}
	
\end{tikzpicture}
\end{center}

\begin{lemma}
\label{lift1} The sets $X_\iota$, $Y_\iota$, and $Z_\iota$ are $\ico \varsigma$-clopen. 
\end{lemma}

\proof
The sets $X_\iota$ are clearly already $\ico 1$-clopen. That the sets $Y_\iota$ are $\ico\varsigma$-clopen follows from the fact that $\varsigma_\iota < \varsigma$; consequently, so too are the sets $Z_\iota$.
\endproof

\begin{lemma}\label{lift2} 
Suppose $\varsigma$ is a limit ordinal. Then, the sets $\{X_\gamma: \gamma < \varsigma\}$ in Definition \ref{cellpartition} form a partition of $e^\varsigma\xi$.
\end{lemma}

\proof 
From \cite[Lemma 2.8]{AguFer15} follows that
\begin{equation*}
\sup_{\iota \in \varsigma \cap \Lim}\alpha_\iota = e^\varsigma\xi,
\end{equation*}
so it suffices to show that 
\begin{equation}\label{toshowlift2}
\lim_{\iota \to \gamma} \alpha_\gamma = e^{\nu\cdot\gamma}(e^\varsigma\zeta + 1) \text{ for each limit } \gamma < \varsigma.
\end{equation}
Write $\gamma = \gamma^* + \omega^\rho$, where $\rho$ is nonzero.  Recall that the functions $e^\iota$ are normal. First, we have:
\begin{align*}
\lim_{\iota \to \gamma} \alpha_\iota = \lim_{\iota \to \gamma} \beta_\iota 
&= \lim_{\iota \to \gamma} e^{\nu\cdot\iota}(\kappa_0 + \kappa_\iota) \\
&\leq \lim_{\iota \to \gamma} e^{\nu\cdot\iota}(\kappa_0 \cdot 2) \\
&\leq \lim_{\iota \to \gamma} e^{\nu\cdot\iota}e(\kappa_0 + 1) \\
&= \lim_{\iota \to \gamma} e^{\nu\cdot\iota}(\kappa_0 + 1)\\
&\leq \lim_{\iota \to \gamma} e^{\nu\cdot\iota}(\kappa_0 + \kappa_\iota),
\end{align*}
so $\lim_{\iota \to \gamma} \alpha_\iota = \lim_{\iota \to \gamma} e^{\nu\cdot\iota}(\kappa_0 + 1)$. 
Using the decomposition of $\gamma$:
\begin{align*}
\lim_{\iota \to \gamma} \alpha_\iota 
&= \lim_{\iota \to \gamma} e^{\nu\cdot\iota}(\kappa_0 + 1)\\
&= \lim_{\iota \to \omega^\rho} e^{\nu\cdot(\gamma^* + \iota)}(\kappa_0 + 1)\\
&= \lim_{\iota \to \omega^\rho} e^{\nu\cdot\gamma^* + \nu\cdot\iota}(\kappa_0 + 1)\\
&= e^{\nu\cdot\gamma^*}\left(\lim_{\iota \to \omega^\rho} e^{\nu\cdot \iota}(\kappa_0 + 1)\right). \numberthis \label{eqPartition0}
\end{align*}
By choice of $\nu$, we have
\[e^\nu(e^\varsigma\zeta + 1) \leq \kappa_0  < e^{\nu + 1}(e^\varsigma\zeta + 1),\]
and thus
\[e^\nu(e^\varsigma\zeta + 1) + 1 \leq \kappa_0 + 1 \leq e^{\nu + 1}(e^\varsigma\zeta + 1),\]
so normality implies that
\begin{align*}
\lim_{\iota \to \omega^\rho} e^{\nu\cdot \iota}\Big(e^\nu (e^\varsigma\zeta + 1) + 1\Big)
&\leq \lim_{\iota \to \omega^\rho} e^{\nu\cdot\iota} (\kappa_0 + 1) \\
&\leq \lim_{\iota \to \omega^\rho} e^{\nu\cdot \iota} e^{\nu +1}(e^\varsigma\zeta + 1)\\
&= \lim_{\iota \to \omega^\rho} e^{\nu\cdot \iota}(e^\varsigma\zeta + 1) \\
&\leq \lim_{\iota \to \omega^\rho} e^{\nu\cdot \iota}\Big(e^\nu (e^\varsigma\zeta + 1) + 1\Big).
\end{align*}
This implies that
\begin{equation}\label{eqPartition1}
\lim_{\iota \to \omega^\rho} e^{\nu\cdot\iota} (\kappa_0 + 1) = \lim_{\iota \to \omega^\rho} e^{\nu\cdot \iota}(e^\varsigma\zeta + 1).
\end{equation}
Now, recall that if $\nu^*$ is a limit, then, for every ordinal $\gamma$, we have
\begin{align*}
e^{\nu^*}(\gamma+1) = \lim_{\iota\to\nu^*}e^\iota(e^{\nu^*}\gamma + 1).
\end{align*}
Using this (for $\nu^* = \nu\cdot\omega^\rho$ and $\gamma = e^\varsigma\zeta$) and the additive indecomposability of $\varsigma$, we obtain:
\begin{align*}
\lim_{\iota \to \omega^\rho} e^{\nu\cdot \iota}(e^\varsigma\zeta + 1)
&= \lim_{\iota \to \nu\cdot \omega^\rho} e^{\iota}(e^\varsigma\zeta + 1)\\
&= \lim_{\iota \to \nu\cdot \omega^\rho} e^{\iota}(e^{\nu\cdot\omega^\rho + \varsigma}\zeta + 1)\\
&= \lim_{\iota \to \nu\cdot \omega^\rho} e^{\iota}(e^{\nu\cdot\omega^\rho} e^\varsigma\zeta + 1)\\
&= e^{\nu\cdot\omega^\rho}(e^\varsigma\zeta + 1).
\end{align*}
Putting this together with equations \eqref{eqPartition0} and \eqref{eqPartition1},
\begin{align*}
e^{\nu\cdot\gamma}(e^\varsigma\zeta + 1) = 
e^{\nu\cdot\gamma^*} e^{\nu\cdot\omega^\rho}(e^\varsigma\zeta + 1)
&= e^{\nu\cdot\gamma^*}\bigg(\lim_{\iota \to \omega^\rho} e^{\nu\cdot \iota}(e^\varsigma\zeta + 1)\bigg)\\
&= e^{\nu\cdot\gamma^*}\bigg(\lim_{\iota \to \omega^\rho} e^{\nu\cdot\iota} (\kappa_0 + 1)\bigg)\\
& = \lim_{\iota \to \gamma} \alpha_\iota ,
\end{align*}
which proves equation \eqref{toshowlift2}. This finishes the proof of the lemma.
\endproof

\subsection{Projections}

\begin{defi}[Projections] \label{defproj}
We define the functions $\pi_0$ and $\pi_1$:

\begin{enumerate}
\item[($\downarrow$)] $\pi_0 \colon [0,e^\varsigma\xi) \to [1, \kappa]$ is defined by:
\begin{equation*}
\pi_0(\alpha) =
\left\{
	\begin{array}{ll}
		1 + \ell^{\varsigma_\iota}\alpha  & \mbox{if } \ell^{\varsigma_\iota}\alpha \leq \kappa_0 \\
		\kappa_1 + \dots + \kappa_{i-1} + 1 + (-(1 + \kappa_0) + \ell^{\varsigma_\iota}\alpha)  & \mbox{otherwise;} 
	\end{array}
\right.
\end{equation*}
where $\iota$ and $i$ are the unique ordinals such that $\alpha \in X_\iota$, $\iota \equiv i \mod m$, and $ i< m$.
\item[($\downarrow$)] The function $\pi_0$ is extended to all of $X_\downarrow$: given $\alpha \in X_\downarrow$, let $\eta$ be least such that $\ell^\eta\alpha<e^\varsigma\xi$. Then, $\pi_0(\alpha) = \pi_0(\ell^\eta\alpha)$.

\item[($\uparrow$)] $\pi_1 \colon X_\uparrow \to [1, \lambda]$ is defined by 
\[\pi_1\alpha = 1 + (-\xi + \ell^\varsigma\alpha).\]
\end{enumerate}
\end{defi}
Clearly, 
\[\pi_1:(X_\uparrow,\ico{\varsigma})\to([1,\lambda],\ico 0)\]
is a surjective d-map. It is not immediately clear whether $\pi_0$ is a d-map; this we verify below.

Observe that, since $X_\downarrow = [1,\Theta]_0 \cap [0, \xi)_\varsigma$ by definition, the ordinal $\eta$ in the second clause above must be strictly smaller than $\varsigma$. 

Let us look at the definition of $\pi_0$ a bit more closely. Consider a typical element $\alpha^*$ of $X_\downarrow$ and generate the sequence
\[\{\ell^\nu\alpha^*:\nu \in Ord\}.\]
Let $\eta$ be least such that $\alpha:= \ell^\eta\alpha^*<e^\varsigma\xi$. Then $\alpha$
 belongs to some cell 
\[X_\iota = [\alpha_\iota, \beta_\iota].\]
By definition, $\beta_\iota$ is of the form $e^{\varsigma_\iota}(\kappa_0 + \kappa_i)$, 
where $\iota \equiv i \mod m$, so
$$\ell^{\varsigma_\iota}\alpha \in [0, \kappa_0 + \kappa_i].$$
Recall that we assumed for simplicity that $\nu \neq 0$ (see p. \pageref{subsectionGLPPartition}).
We distinguished two cases: if $\ell^{\varsigma_\iota}\alpha \in [0, \kappa_0]$, then we had set 
\[\pi_0\alpha^* = \pi_0\alpha = 1 + \ell^{\varsigma_\iota}\alpha \in [1, \kappa_0].\]
Otherwise, $\ell^{\varsigma_\iota}\alpha \in [\kappa_0+1, \kappa_0 + \kappa_i]$, so that
\[(-1 + \kappa_0) + \ell^{\varsigma_\iota}\alpha \in [0,\kappa_i].\]
Its definition then places $\pi_0\alpha^* = \pi_0\alpha$ within the interval
\[[\kappa_1 + \dots + \kappa_{i-1} + 1, \kappa_1 + \dots + \kappa_{i-1} + \kappa_i].\]
Here is the picture:
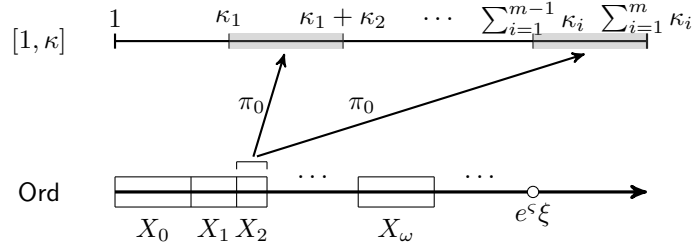
\begin{center}
\begin{tikzpicture}[
    axis/.style={very thick, ->, >=stealth'},
    ll/.style={thick},
    dline/.style={dashed, thin}, 
    pile/.style={thick, ->, >=stealth', shorten <=2pt, shorten
    >=2pt},
    wdot/.style={ 
    circle,
    fill=white,
    draw,
    outer sep=0pt,
    inner sep=1.5pt
  }
    ]
\draw[ll] (0,2) -- (7,2);    
\draw[ll] (0,1.9) -- (0,2.1);  
\draw[ll] (7,1.9) -- (7,2.1);  
\node at (-1,0) {$\Ord$};  
\node at (-1,2) {$[1, \kappa]$};  

\node at (0,2.3) {$1$};  
\node at (1.5,2.3) {$\kappa_1$};  
\draw    (1.5,1.9) -- (1.5,2.1); 
\node at (3,2.3) {$\kappa_1 + \kappa_2$};  
\draw    (3,1.9) -- (3,2.1);
\node at (4.25,2.3) {$\hdots$}; 
\node at (5.5,2.3) {$\sum_{i=1}^{m-1}\kappa_i$};  
\draw    (5.5,1.9) -- (5.5,2.1);
\node at (7,2.3) {$\sum_{i=1}^m\kappa_i$};  
\draw    (7,1.9) -- (7,2.1);
\draw[pile] (1.8,0.4) -- (2.25,1.9);
\node at (1.8,1.15) {$\pi_0$};  
\draw[pile] (1.8,0.45) -- (6.25,1.85);
\node at (3.25,1.15) {$\pi_0$};

\draw (1.6, 0.4) -- (2,0.4);
\draw (1.6, 0.3) -- (1.6,0.4);
\draw (2, 0.3) -- (2,0.4);
    
\draw[axis] (0,0) -- (7,0);

\foreach \x in {0, 1, 1.6, 2}
{
	\draw        (\x,-0.2) -- (\x,0.2);
}
	\node at (2.6,0.2) {$\hdots$};
	
\foreach \x in {3.2, 4.2}
{
	\draw        (\x,-0.2) -- (\x,0.2);
}
	\node at (4.8,0.2) {$\hdots$};
	
	\foreach \x in {-0.2, 0.2}
{
	\draw        (0,\x) -- (2,\x);
	\draw        (3.2,\x) -- (4.2,\x);
}
	\node at (0.5,-0.5) {$X_0$};
	\node at (1.3,-0.5) {$X_1$};
	\node at (1.8,-0.5) {$X_2$};
	\node at (3.7,-0.5) {$X_\omega$};
	
	\node[wdot] at (5.5,0) {};
	\node at (5.5,-0.3) {$e^\varsigma\xi$};
	
	\fill[opacity = 0.3, gray] (1.5,1.9) -- (3,1.9) -- (3,2.1) -- (1.5,2.1) -- (1.5,1.9);
	
	\fill[opacity = 0.3, gray] (5.5,1.9) -- (7,1.9) -- (7,2.1) -- (5.5,2.1) -- (5.5,1.9);
\end{tikzpicture}
\end{center}

\begin{lemma} The projection function $\pi_0$ has the following properties: \label{propproj}
\begin{enumerate}
\item \label{proj1} If $\iota \equiv i \mod m$, then $\pi_0 \colon (X_\iota, \ico \varsigma) \to ([1, \kappa_0] \sqcup [1, \kappa_i], \ico \varsigma)$ is a surjective d-map. If, in addition, $\varsigma = 1$, then it is a homeomorphism.

\item \label{proj2} $\pi_0 \colon (X_\downarrow, \ico \varsigma) \to ([1, \kappa], \ico \varsigma)$ is a surjective d-map.

\item \label{proj3} $\pi_0(\beta_\iota) = \kappa_i$, where $\iota \equiv i \mod m$.
\end{enumerate}
\end{lemma}

\proof
By Lemma \ref{lift1}, the sets $Y_\iota$ and $Z_\iota$ are $\ico\varsigma$-clopen in $X_\iota$. 
Consider generalized cells of the form
\begin{align*}
X_\iota^\eta
&=\big\{
x \in X_\downarrow: 
\ell^{\eta+1}x < e^\varsigma\xi
\text{ and } e^\varsigma\xi \geq \ell^\eta x
\text{ and }
\ell^{\eta + 1}x \in X_\iota
\big\}\\
&=\big\{
x \in X_\downarrow: 
\ell^{\eta+1}x \leq e^\varsigma\xi
\text{ and } e^\varsigma\xi > \ell^\eta x
\text{ and }
\ell^{\eta + 1}x \in X_\iota
\big\}\\
&=
X_\downarrow \cap [0,e^\varsigma\xi]_{\eta+1} \cap (e^\varsigma\xi, \infty)_\eta \cap [\alpha_\iota,\beta_\iota]_{\eta+1}.
\end{align*}
(The second equality follows from the definition of $X_\downarrow$.)
Observe that $[\alpha_\iota,\beta_\iota]_{\eta+1}$ is an $\ico \varsigma$-clopen interval if $\eta < \varsigma$, even when $\alpha_\iota$ is a limit ordinal; this is because $\alpha_\iota$ is always an isolated point in $\ico\varsigma$ by its definition. Similarly, the sets
\begin{align*}
Y^\eta_\iota = X^\eta_\iota \cap [0,\kappa_0]_{\eta + \varsigma_\iota}, Z^\eta_\iota = X^\eta_\iota \cap (\kappa_0,\infty)_{\eta+\varsigma_\iota}
\end{align*}
are $\ico \varsigma$-clopen.
The definition of $\pi_0$ is the same within each $Y^\eta_\iota$ and within each $Z^\eta_\iota$,
and in each of those sets, $\pi_0$ is defined as a combination of additions, substractions, and logarithms and is thus a d-map (recall that $\varsigma$ is additively indecomposable). Additionally, by Lemma \ref{LemmaLeastLog}, if $x \in X_\downarrow$, then the least $\eta$ such that $\ell^\eta x < e^\varsigma\xi$ is a successor ordinal. It follows that the collection of all $Y^\eta_\iota$ and $Z^\eta_\iota$ forms a clopen partition of $X_\downarrow$. Since $\pi_0$ is a d-map on each element of the partition, it is a d-map on all of $X_\downarrow$.

Moreover, if $\varsigma = 1$, then $\cs$ is the constant sequence with value $0$ and so $X_\iota$ is homeomorphic to $[1, \kappa_0] \sqcup [1, \kappa_i]$. This gives item \ref{proj1}. 
Item \ref{proj2} is obtained by a similar argument, as each $X_\iota$ is $\ico\varsigma$-clopen. Item \ref{proj3} follows readily from the definition.
\endproof

\begin{lemma}\label{proj4} 
$\pi_0^{-1}\alpha$ is $\ico \varsigma$-dense in $X_\uparrow$ for any $\alpha \in [1,\kappa]$.
\end{lemma}

\proof
We can even provide witnesses for the density. Let $i$ be least such that 
$$\kappa_1 + \dots + \kappa_{i-1} < \alpha \leq \kappa_1 + \dots + \kappa_i$$
and let $\alpha_0 := -(\kappa_1 + \hdots + \kappa_{i-1} + 1) + \alpha$.

If $\varsigma = 1$, then the result is clear, as any $\ico 1$-neighborhood of any $\beta \in X_\uparrow$ contains an interval of the form $[\delta,\beta]_0$ and, by construction, $\beta$ is a limit of endpoints of cells $X_\iota$. In particular, the interval $[\delta,\beta]$ contains some cell $X_\iota$ with $\iota \equiv i \mod m$ and $\alpha \in \pi_0(X_\iota) = [1, \kappa_0] \sqcup [1, \kappa_i]$ by Lemma \ref{propproj}.\ref{proj1}, from which the result follows.

So suppose $\varsigma \neq 1$. Let $\beta \in X_\uparrow$, so that $\beta$ has $\ico\varsigma$-rank $\rho \geq \xi$ and $U$ be an $\ico\varsigma$-neighborhood of $\beta$.  We distinguish two cases:

\paragraph{Case I:} $\beta = e^\varsigma\xi$; thus $\rho = \xi$. We can apply Lemma \ref{nhbase2} (over the interval topology) to obtain a $\ico\varsigma$-neighborhood base of $e^\varsigma\xi$ consisting of sets of the form
$$(\eta, e^\varsigma\rho]_\gamma$$
for $\eta < e^\varsigma\rho$ and $\gamma < \varsigma$.

Hence, we may assume $U$ is a neighborhood of $e^\varsigma\xi$ of the form $(\eta, e^\varsigma\rho]_\gamma$. We need to find some ordinal $\chi \in U \cap X_\downarrow$ such that $\pi_0\chi = \alpha$. Let $\mu$ be some successor ordinal $\equiv i \mod m$ large enough so that 
\begin{enumerate}
\item $\gamma < \varsigma_\mu < \varsigma$, and
\item $\eta < e^{-\gamma + \varsigma_\mu}(1 + \kappa_0 + \alpha_0)$.
\end{enumerate}
This is certainly possible, as it follows from Lemma \ref{lift2} that:
$$\lim_{\iota \to \varsigma}e^{\varsigma_\iota}(1 + \kappa_0 + \alpha_0) 
= \lim_{\iota \to \varsigma}\alpha_\iota 
= e^\varsigma\xi.$$
\begin{claim}
Let $\chi := e^{\varsigma_\mu}(1 + \kappa_0 + \alpha_0)$. Then $\chi \in U \cap X_\mu$.
\end{claim}
\proof
Clearly, $\chi < e^\varsigma\xi$. By choice of $\mu$, we have
\begin{align} \label{eqmu0}
\eta
< e^{-\gamma + \varsigma_\mu}(1 + \kappa_0 + \alpha_0)
= \ell^\gamma(e^{\varsigma_\mu}(1 + \kappa_0 + \alpha_0)),
\end{align} 
and since $\varsigma_\mu < \varsigma$, we have
\begin{align} \label{eqmu1}
\ell^\gamma e^{\varsigma_\mu}(1 + \kappa_0 + \alpha_0) < e^\varsigma\xi.
\end{align} 
From \eqref{eqmu0} and \eqref{eqmu1} follows that $e^{\varsigma_\mu}(1 + \kappa_0 + \alpha_0) \in U$. It remains to prove that $\chi \in X_\mu$. Denote by $\mu^*$ the immediate predecessor of $\mu$. Notice that
\begin{equation} \label{eqmu}
e^{\varsigma_{\mu^*}}(\kappa_0 + \kappa_{i-1}) 
< e^{\varsigma_{\mu}}(1 + \kappa_0 + \alpha_0) 
\leq e^{\varsigma_\mu}(\kappa_0 + \kappa_i),
\end{equation} 
so that $e^{\varsigma_{\mu}}(1 + \kappa_0 + \alpha_0) \in X_\mu$. This proves the claim.
\endproof
We show that $\pi_0(e^{\varsigma_\mu}(1 + \kappa_0 + \alpha_0)) = \alpha$. Since $\mu \equiv i \mod m$ and 
\begin{equation*}
\kappa_0 < 1 + \kappa_0 + \alpha_0 = 
\ell^{\varsigma_\mu}e^{\varsigma_{\mu}}(1 + \kappa_0 + \alpha_0),
\end{equation*} 
 we have:
\begin{align*}
\pi_0(e^{\varsigma_{\mu}}(1 + \kappa_0 + \alpha_0)) 
& = \kappa_1 + \hdots + \kappa_{i-1} + 1 + (-(1 + \kappa_0) + \ell^{\varsigma_\mu}e^{\varsigma_{\mu}}(1 + \kappa_0 + \alpha_0)) \\
& = \kappa_1 + \hdots + \kappa_{i-1} + 1 + \alpha_0 \\
& = \kappa_1 + \hdots + \kappa_{i-1} + 1 + \left( -(\kappa_1 + \hdots + \kappa_{i-1} + 1) + \alpha \right)\\
& = \alpha.
\end{align*}
Since $U$ was arbitrary, this finishes the proof in this case.

\paragraph{Case II:} $e^\varsigma\xi < \beta$. It is enough to consider the case $\ell^\varsigma\beta = \xi$, as any $\ico\varsigma$-neighborhood of any point of higher rank contains a point of rank $\xi$. 
As observed by Fern\'andez-Duque and Joosten \cite{hyperations}, there is a least $\eta^*<\varsigma$ such that
\[\ell^{\eta^*}\beta = e^\varsigma\xi.\]
Since $\varsigma$ is additively indecomposable, $\eta^*$ must be a successor ordinal, say $\eta+1$. Thus, 
\[\ell^{\eta+1}\beta = e^\varsigma\xi < \ell^\eta\beta,\]
and so $\ell^\eta\beta$ must be of the form $\beta' + e^\varsigma\xi$.
The Hyperexponential Normal Form theorem \index{Hyperexponential Normal Form} \cite[Proposition 2.13]{AguFer15} states that every ordinal $\gamma$ can be uniquely written in the form $e^{\gamma_0}\gamma_1$, where $\gamma_1$ is either additively decomposable or $1$. If $\gamma_1 = 1$, then let us call the expression $e^{\gamma_0}1$ the \emph{normal form expansion} of $\gamma$. Inductively, the normal form expansion of
\[\gamma = e^{\gamma_0}(\gamma_1 + \gamma_2)\]
is defined to be 
\[\texttt{nf}(\gamma) = e^{\gamma_0}(\gamma_1 + \texttt{nf}(\gamma_2)).\]
Observe that $\eta<\varsigma$, for otherwise we would have
\[\xi < e^\varsigma \xi \leq \ell^\varsigma\beta,\]
contradicting the choice of $\beta$.
Thus, if one writes out the normal form expansion of $\beta$, one obtains an expression of the form
\[e^{\beta_0}\Big(\beta_1 + e^{\beta_1}\big(  \dots (\beta' + \texttt{nf}(e^\varsigma\xi)) \dots \big)\Big),\]
where $e^\varsigma\xi < \beta' + e^\varsigma\xi$ and all the exponents to the left of $\beta' +\texttt{nf}(e^\varsigma\xi)$ add up to $\eta$. Consider the sequence $\{\beta(\iota):\iota< e^\varsigma\xi\}$, where $\beta(\iota)$ is the ordinal one obtains if one substitutes $\iota$ for the rightmost occurrence of $\texttt{nf}(e^\varsigma\xi)$ in the normal form expansion of $\beta$ (the occurrence indicated in the equation displayed above). If $U$ is a $\ico\varsigma$-neighborhood of $\beta$, then, by Lemma \ref{nhbase3}, $U$ contains a set of the form 
\[B_r(\beta) = \bigcap_{i \in \text{dom}(r)}(r(i),\ell^i \beta]_i,\]
where $r:\varsigma\to\beta+1$ is a finite partial function.
It follows that every such set $U$, if nonempty, contains cofinally many ordinals of the form $\beta(\iota)$. Since $e^\varsigma\xi< \beta' + e^\varsigma\xi$ and $e^\varsigma\xi$ is additively indecomposable, we must have $e^\varsigma\xi < \beta'$,
so it follows that for each $\beta(\iota)$, 
\[\pi_0\beta(\iota) = \pi_0\ell (\beta' + \iota) = \pi_0 \ell\iota.\]
An argument as in Case I shows that there is some $\iota<e^\varsigma\xi$ such that $\beta(\iota) \in U$ and $\pi_0\beta(\iota) = \alpha$.
\endproof


\subsection{The polytopology} 
It remains to define a $\vec\vartheta$-polytopology $(\mathcal{R}_0,\hdots, \mathcal{R}_n)$ on $([1, \lift], \ico\varsigma)$ such that the projection mappings 
\[\pi_0:(X_\downarrow, \mathcal{R}_i)\to ([1,\kappa],\mathcal{T}_i)\]
and 
\[\pi_1:(X_\uparrow, \mathcal{R}_{i})\to ([1,\lambda],\mathcal{S}_{i})\]
remain d-maps. Recall that $n$ denotes  the length of $\vec\vartheta$. 

To begin, we observe that since 
$$\pi_0: (X_\downarrow, \mathcal{I}_\varsigma) \to ([1, \kappa], \mathcal{I}_\varsigma)$$ 
is a d-map by Lemma \ref{propproj}.\ref{proj2}, we may apply Lemma \ref{Omegapullback} to obtain a $\vec\vartheta$-polytopology $(X_\downarrow, \hat{\mathcal{T}_0}, \hdots, \hat{\mathcal{T}_n})$ over $\ico \varsigma$ such that 
\[\pi_0: (X_\downarrow, \hat{\mathcal{T}_i}) \to ([1,\kappa], \mathcal{T}_i)\] 
is a d-map for each $0\leq i\leq n$.

This $\vec\vartheta$-polytopology $(\hat{\mathcal{T}_0}, \hdots, \hat{\mathcal{T}_n})$ is not, however, a topology on $[1, \lift]$, so we need to extend it.
For each $i$, let $\hat{\mathcal{R}}_i$ be the smallest topology on $[1,\Theta]$ extending $\ico{\varsigma + \vartheta_i}$ and containing all sets in $\hat{\mathcal{T}_i}$. Since $X_\uparrow$ is $\ico{\varsigma+1}$-clopen, $\hat{\mathcal{R}}_i$ is simply equal to $\ico{\varsigma+\vartheta_i}$ when restricted to $X_\uparrow$, for $0<i$. We are closer to our goal, but not done yet, since the space
\[([1,\Theta], \hat{\mathcal{R}}_0,\hdots, \hat{\mathcal{R}}_n)\]
might not be a $\vec\vartheta$-polytopology, as e.g., $\hat{\mathcal{R}}_0$ might not be $\vartheta_1$-maximal around points in $X_\uparrow$.

Note that $-\xi + \ell^\varsigma$ is the rank function of $(X_\uparrow,\ico\varsigma)$ (viewed as a subspace of $([1,\Theta],\ico\varsigma)$) and so
$$-\xi + \ell^\varsigma: (X_\uparrow,\ico\varsigma) \to ([0, -1 +\lambda], \ico0)$$
is a d-map.
Since $[0, -1 + \lambda]$ is homeomorphic to $[1, \lambda]$, it follows that 

$$1 + (- \xi  + \ell^\varsigma): (X_\uparrow,\ico\varsigma) \to ([1, \lambda], \ico0)$$
is also a d-map. But this is precisely equal to $\pi_1$. Having only one point of each rank, the space $([1,\lambda],\ico 0)$ has no proper rank-preserving extensions, and in particular is $\vartheta_1$-maximal. By the claim within the proof of Lemma \ref{Omegapullback}, if $(X_\uparrow, \mathcal{R})$ is any $\vartheta_1$-extension of $(X_\uparrow, \ico\varsigma)$, then 
$$\pi_1: (X_\uparrow,\mathcal{R}) \to ([1, \lambda], \ico0)$$
remains a d-map. Let $([1,\Theta],\mathcal{R}_0)$ be a $\vartheta_1$-maximal extension of $([1,\Theta], \hat{\mathcal{R}_0})$. Then, $\mathcal{R}_0$ only adds neighborhoods around points of  rank some $\rho$ such that $0<\ell^{\vartheta_1}\rho$ and, moreover, only neighborhoods around points in $X_\uparrow$, since $(X_\downarrow, \hat{\mathcal{T}}_0)$ was already $\vartheta_1$-maximal. Given a point $x\in X_\uparrow$, and recalling that $\xi$, the minimum $\ico\varsigma$-rank of points in $X_\uparrow$, is a successor ordinal, we see that
\begin{align*}
0 < \ell^{\vartheta_1}\rho_{([1,\Theta],\mathcal{R}_0)}x 
&\text{ if, and only if, }
0 < \ell^{\vartheta_1}\ell^{\varsigma}x \\
&\text{ if, and only if, }
0 < \ell^{\vartheta_1}\rho_{(X_\uparrow,\ico\varsigma)}x.
\end{align*}
Thus, the space $(X_\uparrow, \mathcal{R}_0)$ is a $\vartheta_1$-extension of $(X_\uparrow,\ico\varsigma)$. It follows that
$$\pi_1: (X_\uparrow,\mathcal{R}_0) \to ([1, \lambda], \ico0)$$
remains a d-map. We may now apply Lemma \ref{Omegapullback} to obtain a $\vec\vartheta$-polytopology 
$$(X_\uparrow, \mathcal{R}_0, \hat{\mathcal{S}_1}, \hdots, \hat{\mathcal{S}_n})$$ 
over $(X_\uparrow,\mathcal{R}_0)$ such that 
\[\pi_1: (X_\uparrow, \hat{\mathcal{S}_i}) \to ([1,\lambda], \mathcal{S}_i)\] 
is a d-map for each $1\leq i\leq n$.
For each $1\leq i\leq n$, we let $\mathcal{R}_i$ be the disjoint union
\[(X_\downarrow, \hat{\mathcal{T}}_i) \sqcup (X_\uparrow, \hat{\mathcal{S}}_i).\]
The sets $X_\uparrow$ and $X_\downarrow$ are $\ic{\vartheta_1}{\mathcal{R}_0}$-clopen and so it follows that
\[([1,\Theta],\mathcal{R}_0,\mathcal{R}_1,\hdots,\mathcal{R}_n)\]
is a $\vec\vartheta$-polytopology. Moreover, we have seen that
\[\pi_0: (X_\downarrow, \hat{\mathcal{T}_i}) \to ([1,\kappa], \mathcal{T}_i)\] 
is a d-map for each $0\leq i\leq n$ and that
\[\pi_1: (X_\uparrow, \mathcal{R}_i) \to ([1,\lambda], \mathcal{S}_i)\] 
is also a d-map for each $0\leq i\leq n$. The other conditions in the statement of the Product Lemma are easy to check from the construction, so its proof is complete.

\bibliographystyle{abbrv}
\bibliography{bibtesis}

\end{document}